\documentclass[11pt,reqno]{amsart} 
\usepackage[margin=1in]{geometry}
\usepackage{enumitem}
\usepackage{amssymb}
\usepackage[linesnumbered,commentsnumbered,ruled,vlined]{algorithm2e}
\usepackage[dvipsnames]{xcolor}
\usepackage{tikz}
\usepackage{tikz-cd}
\usepackage{graphicx}
\usetikzlibrary{patterns,external}
\usetikzlibrary{calc}
\usepackage{natbib}    
\setcitestyle{square,numbers}
\usepackage{hyperref}
\usepackage{cleveref}
\usepackage{verbatim}
\usepackage{blkarray}
\newcommand\scalemath[2]{\scalebox{#1}{\mbox{\ensuremath{\displaystyle #2}}}}

\usetikzlibrary{decorations.markings}
\definecolor{fondo}{rgb}{0.898,0.996,0.898}

\pgfkeys{/tikz/.cd,
  K/.store in=\K,
  K=1   
   }

\usetikzlibrary{arrows.meta,
                quotes}

\title{Lawrence Lifts, Matroids, and Maximum Likelihood Degrees}

\author{Taylor Brysiewicz}
\address{Department of Mathematics\\ University of Western Ontario\\
 2004 Perth Dr, London, ON N6G 2V4, Canada}
\email{tbrysiew@uwo.ca}

\author{Aida Maraj}
\address{Department of Mathematics\\
University of Michigan\\
1855 East Hall, Ann Arbor, MI 48109  \\
 USA}
\email{maraja@umich.edu}

\newcommand{\C}{{\mathbb{C}}}

\newcommand{\Z}{{\mathbb{Z}}}

\newcommand{\p}{\mathbf{p}}
\newcommand{\mydef}[1]{{\color{magenta}#1}}
\newcommand{\mycomment}[1]{}

\usepackage[all,cmtip]{xy}

\newcommand{\x}{\textbf{x}}
\renewcommand{\u}{\textbf{u}}
\renewcommand{\v}{\textbf{v}}
\newcommand{\w}{\textbf{w}}
\newcommand{\s}{\textbf{s}}

\newcommand{\y}{\textbf{y}}
\renewcommand{\r}{\textbf{r}}
\renewcommand{\i}{\textbf{i}}

\renewcommand{\hat}{\widehat}

\theoremstyle{definition}
\newtheorem{theorem}{Theorem}
\newtheorem{definition}[theorem]{Definition}
\newtheorem{lemma}[theorem]{Lemma} 
\newtheorem{corollary}[theorem]{Corollary} 
 
\newtheorem{example}[theorem]{Example}
\newtheorem{remark}[theorem]{Remark}
\newtheorem{proposition}[theorem]{Proposition}

\DeclareMathOperator{\supp}{supp}

\begin{document}
\begin{abstract} We express the  maximum likelihood (ML) degrees of a family toric varieties in terms of M\"obius invariants of matroids. The family of interest are those parametrized by monomial maps given by Lawrence lifts of totally unimodular matrices with even circuits. Specifying these matrices to be vertex-edge incidence matrices of  bipartite graphs gives the ML degrees of  some hierarchical models and  three dimensional quasi-independence models. Included in this list are the no-three-way interaction models with one binary random variable, for which, we give closed formulae. 

\end{abstract}

\maketitle

\vspace{-0.4in}

\section{Introduction}

Given an integer matrix $\mydef{A} \in \mathbb{Z}^{\mydef{d} \times \mydef{n}}$, its \mydef{Lawrence lift} is the $(2d+n) \times 2n$ matrix:
\[
\mydef{\Lambda(A)} = {\footnotesize{ \begin{bmatrix} A & 0 \\ 0 & A \\ I_n & I_n \end{bmatrix}}} \qquad \text{ where } I_n \text{ is the } n\times n \text{ identity matrix}.
\]
The matrix $\Lambda(A)$ encodes a monomial parametrization of the toric variety $\mydef{X_{\Lambda(A)}}$ via its columns. Varieties of this form appear often
\cite{almendra2023markov, Amendola17,bayer2001sysygies, Bernstein15, Bernstein17,  coons2020generalized,de2006markov,diaconis1998algebraic,ishizeki2003standard,santos2003higher,sturmfels1997variation,Sturmfels96, Sullivant.2018,vlach1986conditions}. 
Of particular interest are their \emph{degrees} and \emph{maximum likelihood (ML) degrees}. 
The ML degree of a variety counts the critical points of a monomial function  on it. Although this number originates from the problem of \emph{maximum likelihood estimation} in statistics, 
  the ML degree has found recent popularity beyond statistics, within combinatorics, algebraic geometry, and physics~\cite{LikelihoodDegenerations,Huh,HS14}. Our main result is a combinatorial description of the degree and ML degree of $X_{\Lambda(A)}$, expressed as   evaluations of the Tutte polynomial $\mydef{\tau_{\mathcal M(A)}(a,b)}$ of the column matroid~$\mydef{\mathcal M(A)}$ of the \emph{original} matrix $A$. 
\begin{theorem}
\label{thm:MainTheorem2} Fix an integer matrix $A \in \mathbb{Z}^{d \times n}$. Then,
\begin{align*}
\textrm{deg}(X_{\Lambda(A)}) = \tau_{\mathcal M(A)}(1,1) &\quad \text{ if }A\text{ is a totally unimodular matrix,}\\ 
\textrm{mldeg}(X_{\Lambda(A)}) = \tau_{\mathcal M(A)}(1,0) &\quad  \text{ if }A\text{  is a totally unimodular matrix with only even circuits.}
\end{align*} 
\end{theorem}
The evaluation $\tau_{\mathcal M(A)}(1,1)$ counts the \emph{bases} of $\mathcal M(A)$ whereas $\tau_{\mathcal M(A)}(1,0)$, which is the absolute value of the \emph{M\"obius invariant} of $\mathcal M(A)$,  counts those bases with \emph{external activity zero}. This gives a combinatorial realization of the well-known fact that the degree of a toric variety bounds its ML degree. It also gives a formula for the volume of the convex hull of the columns of $\Lambda(A)$ \cite{Kushnirenko}.

The hypotheses of Theorem \ref{thm:MainTheorem2} are  satisfied, for instance, when $A$ is the vertex-edge incidence matrix $\mydef{A_G}$ of a bipartite graph $\mydef{G}$. There, the bases of $\mathcal M(A_G)$ are spanning forests (i.e. maximal independent subsets of edges). Formulae for the number of spanning forests, as well as those with external activity zero, are known when $G$ is the complete bipartite graph $\mydef{K(m_1,m_2)}$.  The corresponding toric variety $X_{\Lambda(A_{K(m_1,m_2)})}$ corresponds to the \emph{no-three-way interaction model}  $\mydef{M_{2,m_1,m_2}}$. Understanding the  structure, degrees, and  ML degrees of no-three-way interaction models are ongoing challenges in algebraic statistics \cite{almendra2023markov, de2006markov, santos2003higher, Amendola17, coons2020generalized}.  Specializing \Cref{thm:MainTheorem2} to this setting answers conjectures in \cite[Conjecture 36]{Amendola17} and \cite[Section 7]{coons2020generalized} affirmatively.

\begin{corollary}
\label{cor:MainCorollary1}
Let $m_1,m_2 \in \mathbb{N}$ and $G=K(m_1,m_2)$ be the complete bipartite graph. Then
\[
\deg(X_{\Lambda(A_{G})}) = m_1^{m_2-1}m_2^{m_1-1} \quad \quad  \text{and}  \quad  \quad 
\textrm{mldeg}(X_{\Lambda(A_{G})}) = \sum_{k=1}^{m_1} \left( \frac{1}{k} \sum_{i=1}^k (-1)^{m_1-i}{{k}\choose{i}}i^{m_1} \right)k^{m_2}.
\]
In particular, these are the degree and ML degree of the no-three-way interaction model $M_{2,m_1,m_2}$.
\end{corollary}

\noindent \textbf{Structure of the paper:} We first prove the result for the degree of $X_{\Lambda(A)}$ in Section \ref{sec:degree}. Prior to our proof, we give background on toric varieties and matroids. Since $A$ is totally unimodular, its circuits index a Gr\"obner basis of binomials with squarefree terms. Its initial ideal is a squarefree monomial ideal whose degree is determined geometrically by counting the coordinate linear spaces it cuts out. 

In Section \ref{sec:mldegree}, we turn our attention to the ML degree of $X_{\Lambda(A)}$. Although our main result requires total unimodularity of $A$, several lemmas we prove along the way do not. We show the solutions to the maximum likelihood equations for $X_{\Lambda(A)}$ are supported on the torus, even when $A$ is not totally unimodular (Lemma~\ref{lem:SupportedOnTorus}). Moreover, the Lawrence structure present in the polynomial system strongly suggests a way to eliminate half the variables (Corollary~\ref{cor:mldegEliminatedNEW}). What remains are two sets of equations: one set which cuts out a scheme $Y$ and another which defines a linear space $L$ of complementary dimension to $Y$. When $A$ is totally unimodular with even circuits, we identify the part of $Y$ \emph{at infinity} with some reciprocal linear space $\mathcal L^{-1}$ whose underlying matroid is $\mathcal M(A)$  (Lemma~\ref{lem:reciprocalRecognitionNEW}). We use the main result of \cite{ProudfootSpeyer} to show that the degree of $\mathcal L^{-1}$ is $\tau_{\mathcal M(A)}(1,0)$ (\Cref{lem:proudfootspeyerNEW}). Showing that $\mathcal L^{-1}$ does not intersect the part of $L$ at infinity (Theorem~\ref{thm:nothingAtInfinityNEW}) implies that $L$ is sufficiently generic so that $|L \cap Y|$ is the degree of $Y$, finishing the proof.  We end \Cref{sec:mldegree} with a  characterization of totally unimodular matrices $A$ with even circuits for which the ML degree of $X_{\Lambda(A)}$ is one (\Cref{prop:MLdeg1}), and hence the maximum likelihood estimator is a rational function in the sufficient statistics. 

In Section \ref{sec:graphs} we place our results in the context of graph theory, where $A$ is an incidence matrix~$A_G$ of a graph $G$. We show that our main theorem applies when $G$ is bipartite (\Cref{thm:mld_graphs}). We extend this interpretation to directed and signed graphs. 

We discuss the application of our results to algebraic statistics in Section \ref{sec:statistics}. We explain how toric varieties appear in algebraic statistics as log-linear models. We focus on those log-linear models for which our main result applies. This includes no-three-way interaction models, three dimensional quasi-independence models, and hierarchical models. 
 We end by connecting our results to those regarding the ML degrees of Gaussian models with linear-diagonal covariance matrices, and linear models (\Cref{rem:Gaussian}); their ML degrees are all (absolute values of) M\"obius invariants of matroids.

\section{Toric Varieties and Degrees}
\label{sec:degree}
Fix a matrix  $\mydef{A} \in \mathbb{Z}^{\mydef{d} \times \mydef{n}}$ with columns $\mydef{a_1},\ldots,\mydef{a_n}$ with associated  monomial map:
\begin{align*}
\mydef{\varphi_A}: \mathbb{C}[x_1,\ldots,x_n] &\to \mathbb{C}[t_1^{\pm},\ldots,t_d^{\pm}] \\
x_j &\mapsto \prod_{i=1}^d t_i^{a_{ij}}.
\end{align*} 
Define \mydef{$I_{A}$}$\subseteq \mydef{\C[\x]}=\mathbb{C}[x_1,\ldots,x_n]$ to be the kernel of $\varphi_A$. The zero-set of the ideal $I_A$, considered in the algebraic torus 
\[\mydef{(\mathbb{C}^\times)^n}=\{\mydef{\x}=(x_1,\ldots,x_n) \in \mathbb{C}^n \mid x_i \neq 0 \text{ for all }i=1,\ldots,n\},
\] is an \mydef{(affine) toric variety}:
\begin{equation}
\label{eq:toricvariety}
\mydef{X_A} = {\mydef{\mathcal V^\times(I_A)}} = \{\x\in (\mathbb{C}^\times)^n \mid f(\textbf{x}) = 0 \text{ for all } f \in I_{A}\} \subseteq (\mathbb{C}^\times)^n.
\end{equation}
There is an intimate relationship between properties of $X_A$ and those of $A$; for example, the dimension of $X_A$ is the rank of $A$. For background on toric varieties, see \cite{Cox11}.
 \begin{remark}
 We will use $\mydef{\mathcal V}$ instead of $\mathcal V^\times$ when considering the ambient spaces of  $\mathbb{C}^n$ or $\mathbb{P}^{n-1}$.
 \end{remark}

We recall some fundamental facts about $X_A$ and its equations from \cite[Section 4]{Sturmfels96}.
The \mydef{support} of a vector $\mydef{\v}=(v_1,\ldots,v_n) \in \mathbb{C}^n$ is the subset $\mydef{\textrm{supp}(\v)} = \{i \mid v_i \neq 0\} \subseteq \{1,2,\ldots,n\}$. An integer vector $\v$ in the kernel of $A$ is a \mydef{circuit} of $A$ if (i) its coordinates are relatively prime and (ii) within the kernel, $\v$ has minimal support. We write the set of all circuits of $A$ as $\mydef{\mathcal C_A}$.  We say a circuit~$\v$ is \mydef{even/odd} if $\mydef{|\v|}=|\v|_1=|\supp(\v)|$ is even/odd respectively. We write integer vectors as the difference of two nonnegative vectors like $\v = \mydef{\v_+}-\mydef{\v_-}$. Given a circuit $\v \in \mathcal C_A$ we define the binomial $\mydef{f_{\v}(\x)}=\x^{\v_+}-\x^{\v_-} \in \mathbb{C}[\x]$.

The following proposition shows that polynomials indexed by the circuits $\mathcal C_A$ generate  $I_A$. 
\begin{proposition}\cite[Cor 4.3]{Sturmfels96}
\label{prop:Binomials} For any $A \in \mathbb{Z}^{d \times n}$ the ideal $I_A$ is generated by binomials:
\[
I_A = \langle f_\v(\x) \mid \v \in \mathcal C_A \rangle.
\]
\end{proposition}

We define the degree of a variety/scheme geometrically. The \mydef{degree} of an irreducible variety $X$ is the  number of points in an intersection of $X$ with a generic linear space of complementary dimension. When $X$ has more than one component, its degree is the sum of the degrees of the top-dimensional components. When $X$ is a nonreduced scheme, its degree is the sum of the degrees of its top-dimensional components, counted with the appropriate multiplicities. 
On the algebraic side of things, we define the degree of an ideal to be the degree of the scheme it cuts out. Algorithmically, one may compute this degree using \emph{Gr\"obner bases} as explained below. 

Recall that a subset $\mathcal G = \{g_1,\ldots,g_k\}$ of an ideal $I$ is a \mydef{Gr\"obner basis} with respect to the monomial order $\prec$ if the \mydef{initial ideal} $\mydef{\textrm{in}_{\prec}(I)} = \langle \textrm{in}_{\prec}(f) \mid f \in I\rangle $
is equal to the ideal $\langle \textrm{in}_{\prec}(g_i) \mid i=1,\ldots,k  \rangle$ generated by initial terms of $\mathcal G$. A \mydef{universal Gr\"obner basis} is a Gr\"obner basis with respect to every ordering. A monomial order $\prec$ is a \mydef{degree ordering} if it refines the partial order of \emph{total degree} on monomials. If $\mathcal G$ is a Gr\"obner basis of $I$ with respect to a degree ordering $\prec$, the {degree} of $I$ is the degree of the initial ideal $\textrm{in}_{\prec}(I)$.

An integer matrix $A$ is called \mydef{totally unimodular} when each of its minors are either $0$, $1$, or $-1$. When $A$ is totally unimodular, the generators of $X_{A}$ in Proposition \ref{prop:Binomials} enjoy additional properties. 
\begin{proposition}\cite[Prop 4.11, and Prop 8.11]{Sturmfels96}
\label{prop:Sturmfels2} If $A \in \mathbb{Z}^{d \times n}$ is totally unimodular, then the binomials $\{f_{\v}(\x)\}_{\v \in \mathcal C_A}$ form a universal Gr\"obner basis for $I_A$. Moreover, the monomials involved in each binomial $f_{\v}(\x)$ are squarefree.
\end{proposition}

\Cref{prop:Sturmfels2} opens the door for a combinatorial interpretation of the {degree} of $X_A$ when $A$ is totally unimodular: it is equal to the degree of a squarefree monomial ideal.  Squarefree monomial ideals are intersections of monomial prime ideals, and so the schemes they cut out are reduced unions of \mydef{coordinate linear spaces}, that is, those of the form 
\[
\mydef{L_S} = \{\x \in \mathbb{C}^n \mid x_i \neq 0 \text{ for all } i \in S\} \quad \quad \text{ for }\quad \quad  S \subseteq \{1,2,\ldots,n\}.
\]

\begin{lemma}
\label{lem:combinatorialDegree}
Suppose $\mathcal G=\{g_1,\ldots,g_k\}\subseteq \mathbb{C}[\x]$ is a Gr\"obner basis with respect to a degree ordering $\prec$  and each initial monomial $\x^{\alpha_i}$ is squarefree. Then the degree of the ideal generated by $\mathcal G$ is the number of maximal subsets of $\{1,2,\ldots,n\}$ which do not contain any subset in $\{\textrm{supp}(\alpha_i)\}_{i=1}^k$.
\end{lemma}
\begin{proof}
Since $\mathcal G$ is a Gr\"obner basis with respect to a term order that respects degree, the degree of the ideal $I$ generated by $\mathcal G$ is equal to the degree of the monomial ideal $\langle \x^{\alpha_i} \mid i=1,\ldots, k\rangle$. The variety it cuts out is a union of coordinate linear spaces $L_{S_1}\cup\cdots \cup L_{S_r}$ and so its degree is the number of such linear spaces of maximal dimension. Equivalently, this is the number of subsets~$S \subseteq \{1,2,\ldots,n\}$ which are maximal amongst those which contain no element of  $\{\textrm{supp}(\alpha_i)\}_{i=1}^k$. 
\end{proof}

We turn to our main object of interest: the toric variety $X_{\Lambda(A)}$ coming from the \mydef{Lawrence lift} 
\[
\mydef{\Lambda(A)} = \begin{bmatrix} A & 0 \\ 0 & A \\ I_n & I_n\end{bmatrix}.
\] When working with the toric variety $X_{\Lambda(A)}\subseteq \mathbb{C}^{2n}$ we double-up our variables and use the coordinates $(x_1,\ldots,x_n,y_1,\ldots,y_n) = (\x,\mydef{\y})$.  
The circuits of $\Lambda(A)$ are in bijection with the circuits of $A$: indeed, $\v \in \mathcal C_{A}$ if and only if $(\v,-\v) \in \mathcal C_{\Lambda(A)}$. We make use of the shorthand notation
\[
\mydef{f_{\v}(\x,\y)} = f_{(\v,-\v)}(\x,\y) =\x^{\v_+}\y^{\v_-}-\x^{\v_-}\y^{\v_+} \text{ for } \v \in \mathcal C_A.
\]
We remark that  $\Lambda(A)$ is totally unimodular whenever $A$ is. Consequently, for totally unimodular~$A$, Proposition \ref{prop:Sturmfels2} applies to $\Lambda(A)$ and so the hypotheses of Lemma \ref{lem:combinatorialDegree} are satisfied for the Gr\"obner basis~$\{f_{\v}(\x,\y)\}_{\v \in \mathcal C_A}$ of $I_{\Lambda(A)}$. Thus, to determine $\textrm{deg}(X_{\Lambda(A)})$, it is enough to count the number of subsets of columns of $\Lambda(A)$ which do not contain lead monomials of that Gr\"obner basis.

\noindent\textbf{Matroids:} The supports of the circuits $\mathcal C_A$ define the circuits $\mathcal C_{\mathcal M(A)}$ of a matroid $\mydef{\mathcal M(A)}$ on the ground set $\mydef{E} = \{1,2,\ldots,n\}$.  The \mydef{independent} sets in $\mathcal M(A)$ are those which do not contain a circuit. The \mydef{bases} of $\mathcal M(A)$ are the maximal independent subsets of $E$. The \mydef{rank} of a subset $S \subseteq E$ is the maximum size of an independent subset of $S$. 
Associated to the matroid $\mathcal M(A)$, is a bivariate polynomial 
\[
\mydef{\tau_{\mathcal M(A)}(a,b)} = \sum_{S \subseteq E} (a-1)^{\textrm{rank}(\mathcal M(A))-\textrm{rank}(S)}(b-1)^{|S|-\textrm{rank}(S)}= \sum_{i,j} \mydef{t_{ij}}a^ib^j,
\] called the \mydef{Tutte polynomial} of $\mathcal M(A)$. The evaluation $\tau_{\mathcal M(A)}(1,1)$ is the number of bases of $\mathcal M(A)$. Bases of a matroid may be partitioned depending on their ``internal and external activities''. To define these terms, we first fix a total order $\leq$ on the ground set $E$.

Fix a basis $B$ of $\mathcal M(A)$.  An element $e\in B$ is \mydef{internally active} with respect to $B$ if $e$ is the smallest element of its fundamental co-circuit. More importantly for our purposes is the definition of \emph{external activity}. An element $e \in E-B$ is \mydef{externally active} in the basis $B$ if $e$ is the smallest element of its fundamental circuit with respect to $B$, i.e. the unique circuit contained in $B\cup\{e\}$. The coefficients $t_{ij}$ in the definition of the Tutte polynomial count the number of bases with internal activity $i$ and external activity $j$. See \Cref{ex:externallyzerospanning} and \Cref{fig:spanning_trees} for explicit examples in the case that $\mathcal M(A)$ is a matroid coming from a graph.

 Similar to how we doubled the variables for $X_{\Lambda(A)}$, we index the columns of $\Lambda(A)$ by the set $\mydef{E_{\x,\y}}=\{\mydef{1_\x},\ldots,\mydef{n_\x},\mydef{1_\y},\ldots,\mydef{n_\y}\}$. This is the ground set of the matroid $\mathcal M(\Lambda(A))$. We refer to the first $n$ elements of $E_{\x,\y}$ as $\mydef{E_\x}$ and the last $n$ as $\mydef{E_\y}$.

\begin{theorem}
\label{thm:degreeLawrence}
Let  $A \in \mathbb{Z}^{d \times n}$ be a totally unimodular matrix. Fix a monomial order $\prec$.
The number of maximal subsets of $E_{\x,\y}$ which do not contain any initial term of $\{f_{\v}(\x,\y)\}_{\v \in \mathcal C_A}$ is equal to the number of bases of $\mathcal M(A)$. Consequently,
\[
\deg(X_{\Lambda(A)}) = \tau_{\mathcal M(A)}(1,1).
\]
\end{theorem}
\begin{proof}
We construct such maximal subsets explicitly. Let $\prec$ be any monomial order which respects degree. Let $B$ be any basis of $\mathcal M(A)$. Given $S \subseteq \{1,2,\ldots,n\}$ we write $S_\x$ and $S_\y$ for the corresponding subsets of $E_{\x,\y}$. Since $B$ contains no circuit of $\mathcal M(A)$, the set $B_{\x,\y}=B_{\x} \cup B_{\y}$ contains the support of no monomial appearing in any binomial $f_{\v}(\x,\y)$ for $\v \in \mathcal C_A$. In particular, it does not contain the support of any initial monomial.  For each of the remaining elements $e_\x$ or $e_\y$ of $E_{\x,\y}\setminus B_{\x,\y}$, the set $e \cup B$ contains a unique circuit $\gamma$ of $\mathcal M(A)$ corresponding to a circuit $\v \in \mathcal C_A$. Thus exactly one of $e_\x \cup B_{\x,\y}$ or $e_\y \cup B_{\x,\y}$ contains the initial monomial of $f_\v(\x,\y)$. Appending to~$B_{\x,\y}$ the element of $\{e_\x,e_\y\}$ which does \emph{not}, for each $e \in E\setminus B$, produces a maximal subset $T_B$ of~$E_{\x,\y}$ which contains the support of no initial monomial. In the other direction, given $T_B$, the basis $B$ of $\mathcal M(A)$ may be recovered by taking those $e \in E$ such that $e_\x$ and $e_\y$ are both in $T_B$. 
Applying Lemma \ref{lem:combinatorialDegree} completes the proof.
\end{proof}

\begin{remark} \Cref{thm:degreeLawrence} may also be obtained by combining Corollaries 3.6 and 3.8 of \cite{Dall} in the case that $A$ has full rank. The approach there is through the volume of the convex hull of the columns of $\Lambda(A)$. 
\end{remark}

\begin{example}
Consider the totally unimodular matrix 
\[
A = \begin{bmatrix} -1 & 0 & 0 & 1 & 1 \\ 1 & 1 & 0 & 0 & 0 \\ 0 & -1 & 1 & 0 & -1 \\ 0 & 0 & -1 & -1 & 0  \end{bmatrix} \quad \quad \mathcal C_A =\pm \{(-1,1,1,-1,0),(1,-1,0,0,1),(0,0,1,-1,1)\}.
\]
In $\mathbb{C}[x_1,\ldots,x_5,y_1,\ldots,y_5]$ the toric variety $X_{\Lambda(A)}$ is cut out by the equations
\begin{align*}
f_{(-1,1,1,-1,0)}(\x,\y) = x_2x_3y_1y_4-x_1x_4y_2y_3, \quad &\quad \quad \quad \quad 
f_{(1,-1,0,0,1)}(\x,\y) &= x_1x_5y_2-x_2y_1y_5, \\
f_{(0,0,1,-1,1)}(\x,\y) &= x_3x_5y_4-x_4y_3y_5.
\end{align*}
We remark that there are eight bases of $\mathcal M(A)$:
\[
\text{Bases of }\mathcal M(A):\,\, \{123,124,134,135,145,234,235,245\}.
\]
To see the pheonomenon described in \Cref{thm:degreeLawrence}, we decompose a lead ideal with respect to the degree lexicographic order and obtain
\[
\langle x_1x_4y_2y_3, x_1x_5y_2, x_3x_5y_4 \rangle = 
\langle x_1,x_3 \rangle \cap
\langle x_1,x_5 \rangle \cap
\langle x_4,x_5 \rangle \cap
\langle x_3,y_2 \rangle \cap
\langle x_5,y_2 \rangle \cap
\langle x_5,y_3 \rangle \cap
\langle x_1,y_4 \rangle \cap
\langle y_2,y_4 \rangle.
\]
We see directly that $\textrm{deg}(X_{\Lambda(A)})=8$. 
Each of these prime ideals is obtained by first picking a basis $B$ of $\mathcal M(A)$, starting with all variables indexed by $B$, and appending to that set $x_i$ or $y_i$, depending on which does \emph{not} appear in the lead term of the polynomial $f_{\v}$ indexed by the fundamental circuit in $B \cup \{i\}$. Taking the ideal generated by the complement of the constructed set of variables gives a coordinate ideal.

For example, take the basis $B=\{1,2,3\}$ and start with $T_B=\{x_1,x_2,x_3,y_1,y_2,y_3\}$. 
We consider the variables $x_4,y_4$ indexed by $4 \in \{1,2,\ldots,5\}=E$. The unique circuit in $B\cup\{4\}$ is $\{1,2,3,4\}$, corresponding to the binomial $x_2x_3y_1y_4-x_1x_4y_2y_3$ whose lead monomial is $x_1x_4y_2y_3$. We see that including $y_4$ in $T_B$ preserves the property that $T_B$ does not contain the support of any lead monomial. Similarly, for $5 \in E$, we append $y_5$ to $T_B$, giving $T_B = \{x_1,x_2,x_3,y_1,y_2,y_3,y_4,y_5\}$ whose complement is $\{x_4,x_5\}$. We see that $\langle x_4,x_5\rangle$ is one of the eight prime ideals listed above. 
\end{example}

\section{Maximum Likelihood Degrees}
\label{sec:mldegree}The ML degree  counts the number of smooth critical points of a monomial function on a variety in the torus. We discuss the statistical point-of-view in Section \ref{sec:statistics}. For now, we use the formulation from 
\cite[Corollary~7.3.9]{Sullivant.2018} as the generic number of solutions to a polynomial system.

\begin{definition}
\label{def:mlequationsNEW}
The \mydef{maximum likelihood degree} of a toric variety $X_A = \mathcal V^\times(I_A) \subseteq (\mathbb{C}^\times)^n$ is the number $\mydef{\textrm{mldeg}(X_A)}$ of points in $X_{A}$, counted with multiplicity,  which satisfy the linear equations
\begin{equation}
\label{eq:mlequationsNEW}
A\x=A\u 
\end{equation}
and the inequation $\sum\limits_{i=1}^n x_i \neq 0$ for generic  $\u \in \mathbb{C}^n$. 
\end{definition}

Implicit in Definition \ref{def:mlequationsNEW} is that only solutions in $(\mathbb{C}^\times)^n$ whose coordinates do not sum to zero are counted.  The polynomial equations involved are those which cut out $X_A$ (see Proposition \ref{prop:Binomials}) along with the linear equations $A\x = A\u$. We call these equations, together,  the \mydef{maximum likelihood equations}.
We simplify Definition \ref{def:mlequationsNEW} when the toric variety of interest comes from a Lawrence lift. 
\begin{theorem}
\label{thm:mldNoInequationsNEW}
Let $A \in \mathbb{Z}^{d \times n}$ be an integer matrix and $\Lambda(A) \in \mathbb{Z}^{(2d+n) \times 2n}$ its Lawrence lift. Then $\textrm{mldeg}(X_{\Lambda(A)})$ is the number of solutions in $\mathbb{C}^{2n}$, counted with multiplicity, to the system
\begin{equation}
\label{eq:mldegNoInequationsNEW}
\{f_{\v}(\x,\y) = 0\}_{\v \in \mathcal C_{A}}, \quad \quad A\x = A\u, \quad \quad \x+\y = \u+\w,
\end{equation}
for generic $\u,\w \in \mathbb{C}^n$.
\end{theorem}
This theorem is more than just a combination of Proposition \ref{prop:Binomials} and Definition \ref{def:mlequationsNEW}. Observe that for generic $\u,\w \in \mathbb{C}^n$, \Cref{thm:mldNoInequationsNEW} \emph{implies} that the inequation of Definition \ref{def:mlequationsNEW} is  satisfied and that all solutions occur within the torus $(\mathbb{C}^\times)^{2n}$. 
We show these two facts in  \Cref{lem:sumnotzero} and 
\Cref{lem:SupportedOnTorus} respectively. 

We first establish some notation by writing the linear equations of Definition \ref{eq:mlequationsNEW} for $\Lambda(A)$ as 
\[\Lambda(A) \begin{pmatrix} \x \\ \y \end{pmatrix} = \Lambda(A) \begin{pmatrix} \u \\ \w \end{pmatrix}.\]
There are three parts to these linear equations:
\begin{equation}
\label{eq:mldegLinearSplitNEW}
A\x = A\u, \quad \quad \quad A\y = A\w, \quad \quad \quad \x+\y=\u+\w.
\end{equation}
Altogether, these equations are dependent: the last set along with either of the first two imply the remaining set. In fact, the last set of equations immediately implies that the  inequation of Definition \ref{def:mlequationsNEW} is generically satisfied, as stated below.

\begin{lemma}
\label{lem:sumnotzero}
For generic $\u,\w \in \mathbb{C}^n$, the solutions to \eqref{eq:mldegNoInequationsNEW}  satisfy $\sum\limits_{i=1}^n (x_i+y_i)\neq 0$.
\end{lemma}
\begin{proof}
Pick generic $\u,\w \in \mathbb{C}^n$ and observe that $\sum\limits_{i=1}^n (x_i+y_i) = \sum\limits_{i=1}^n (u_i+w_i) \neq 0$.
\end{proof}
We rely on an elementary linear algebra fact to show the other half of \Cref{thm:mldNoInequationsNEW}.
\begin{lemma}
\label{lem:linearalgebrafactNEW}
Let $A_1 \in \textrm{Mat}_{d,n_1}(\C)$ and $A_2 \in \textrm{Mat}_{d,n_2}(\C)$. If $\textrm{rank}(A_1) < \textrm{rank}(A_2)$ then for any $\textbf{b} \in \mathbb{C}^d$ the linear system~$A_1\x=A_2\u-\textbf{b}$ is inconsistent when~$\u\in \mathbb{C}^n$ is chosen generically. 
\end{lemma}
\begin{proof}
The image of $\x \mapsto A_1\x$ is a $\textrm{rank}(A_1)$-plane $P_1$ and the image of $\u \mapsto A_2\u-\textbf{b}$ is an affine $\textrm{rank}(A_2)$-plane, $P_2$. If $\textrm{rank}(A_1) < \textrm{rank}(A_2)$ then a generic point in $P_2$ is not in $P_1$.
\end{proof}

\begin{lemma}
\label{lem:SupportedOnTorus}
For generic $\u,\w \in \mathbb{C}^n$ all solutions to \eqref{eq:mldegNoInequationsNEW} are contained in $(\mathbb{C}^\times)^{2n}$. 
\end{lemma}
\begin{proof}
Suppose $(\x,\y)$ is a solution to \eqref{eq:mldegNoInequationsNEW}. Based on  $\x$, partition $\{1,2,\ldots,n\}$ into two subsets:
\[
\mydef{S} = \{e \in \{1,2,\ldots,n\} \mid x_e \neq 0\} \quad \quad \quad \mydef{Z}=\{e \in \{1,2,\ldots,n\} \mid x_e = 0\} = \{1,2,\ldots,n\}-S.
\]
Similarly, we partition $\{1,2,\ldots,n\}$ based on $\y$: 
\[
\mydef{S'} = \{e \in \{1,2,\ldots,n\} \mid y_e \neq 0\} \quad \quad \quad \mydef{Z'}=\{e \in \{1,2,\ldots,n\} \mid y_e = 0\} = \{1,2,\ldots,n\}-S.
\]
We now consider the linear equations $A\x = A\u$. Note that $x_e=u_e+w_e$ for all $e \in Z'$, and $A\x = A_S\x_S$ where $\mydef{A_S}$ is the submatrix of $A$ whose columns are indexed by $S$. We rewrite  
\[A\x = A\u \quad \quad \implies\quad \quad  A_S\x_{S}=A\u\quad \quad  \implies\quad \quad  
A_{S-Z'}\x_{S-Z'}=A\u - \underbrace{A_{S \cap Z'} \x_{S \cap Z'}}_{\textbf{b}}.
\vspace{-0.1in}
\]
Crucially, $\textbf{b}$ may be taken to be any vector in $\mathbb{R}^d$ since the vector $\w \in \mathbb{C}^n$ is independent from $\u$. We conclude by Lemma \ref{lem:linearalgebrafactNEW} that $A_{S-Z'}$ must have the same rank as $A$. Let $B \subseteq \{1,2,\ldots,n\}$ be a basis of $\mathcal M(A_{S-Z'})$ which witnesses this rank equality. 

Now suppose towards a contradiction that $i \in Z$. Let $\gamma$ be the unique circuit of $\mathcal M(A)$ contained in $B \cup \{i\}$. This matroid circuit corresponds to two circuits $(-\v,\v),(\v,-\v) \in \mathcal C_{\Lambda(A)}$ for $\v \in \mathcal C_{A}$. Choose the circuit for which $i \in \v_-$ and consider $f_{\v}(\x,\y)$. Since $x_i=0$ we have that 
\[
f_{\v}(\x,\y) = \x^{\v_+}\y^{\v_-}-\x^{\v_-}\y^{\v_+} = 0 \quad \iff \quad  \x^{\v_+}\y^{\v_-}=0.
\]
There are two ways for this to be zero: either $x_j=0$ for some $j \in \textrm{supp}(\v_+)$ or $y_j=0$ for some~$j \in \textrm{supp}(\v_-)$. The former case is impossible since \[\textrm{supp}(\v_+) \subseteq (B \cup \{i\})-\{i\} \subseteq S = \{e \in \{1,2,\ldots,n\} \mid x_e \neq 0\}.\]
But the latter case is impossible, too,  since
\[
\textrm{supp}(\v_-) \subseteq B \subseteq S-Z' \subseteq S' = \{e \in \{1,2,\ldots,n\} \mid y_e \neq 0\}.
\]
Hence, $Z$ must be empty.
\end{proof}

The last equation in \eqref{eq:mldegNoInequationsNEW} suggests a simple elimination via $\y \mapsto \u+\w-\x$. For $\v \in \mathcal C_{A}$ we define 
\[\mydef{g_{\v}(\x,\u,\w)} = f_{\v}(\x,\u+\w-\x),\]
and obtain the following corollary of Theorem \ref{thm:mldNoInequationsNEW}.

\begin{corollary}
\label{cor:mldegEliminatedNEW}
Let $A$ be a $d \times n$ matrix and $\Lambda(A)$ its Lawrence lift. The maximum likelihood degree of $X_{\Lambda(A)}$ is the number of solutions in $\mathbb{C}^n$, counted with multiplicity, to the system
\begin{equation}
\label{eq:mldegEliminatedNEW}
\{g_{\v}(\x,\u,\w)=0\}_{\v \in \mathcal C_A} \quad \quad A\x=A\u
\end{equation}
for generic $\u,\w \in \mathbb{C}^n$.
\end{corollary}
Corollary \ref{cor:mldegEliminatedNEW} now recognizes the maximum likelihood degree of $X_{\Lambda(A)}$ as the degree of a zero-dimensional linear intersection of the scheme $\mydef{Y_A(\u,\w)}$ cut out by $\{g_{\v}(\u,\w)\}_{\v \in \mathcal C_A}$. Unfortunately, these generators do not necessarily form a Gr\"obner basis, and so analyzing the lead ideal may not uncover its degree. 

Instead, we identify the part of $Y_A(\u,\w)$ \emph{at infinity}, and ascertain that degree. Writing $\mydef{\overline{Y_A(\u,\w)}} \subseteq \mathbb{P}^n$ for the projective closure of $Y_A(\u,\w)$ under the inclusion $\x \mapsto [z:x_1:\cdots:x_n]$ of $\mathbb{C}^n \hookrightarrow \mathbb{P}^n$, we call~$\mydef{\mathcal H_\infty}=\mathcal V(z)$ the \mydef{hyperplane at infinity} and $\mydef{Y^{\infty}_A(\u,\w)}=\overline{Y_A(\u,\w)} \cap \mathcal H_\infty$ the \mydef{part of $Y_A(\u,\w)$ at infinity}. We relate this setup to the scheme $\mydef{\widehat{Y}_A(\u,\w)}$  cut out by the homogenizations
$
\{\mydef{\widehat{g}_{\v}(\u,\w)}\}_{\v \in \mathcal C_A}
$ with respect to the variable $z$. 

One may think about the ideal generated by $\{ \widehat{g}_{\v}(\u,\w)  \}_{ \v \in \mathcal C_A}$ as a \mydef{n\"aive homogenization} and~$\widehat{Y}_A(\u,\w)$ as a \mydef{n\"aive closure}. Indeed, these need not agree with $\overline{Y_{A}(\u,\w)}$ or the ideal of polynomials vanishing on it; but there is a relationship. The scheme $\widehat{Y}_A(\u,\w)$ has every component of $\overline{Y_A(\u,\w)}$ as its own component, however, the n\"aive closure may have additional components supported on $\mathcal H_\infty$ \cite[Lemma 2.1]{Stiefel}. If the n\"aive homogenization is particularly well-behaved at infinity, we can conclude that its part at infinity coincides with the part of $\overline{Y_A(\u,\w)}$ at infinity.

\begin{proposition}
\label{prop:NaiveHomogenizationDegree}
If the ideal $J=\langle \widehat{g}_{\v}(\u,\w) \rangle+\langle z \rangle$ is radical of dimension $\dim(\overline{Y_A(\u,\w)})-1$, then~$\overline{Y_A(\u,\w)}$ has a unique component of top dimension, and this component is reduced. Consequently, the part of $\overline{Y_A(\u,\w)}$ at infinity is cut out by $J$ and the degree of $J$ is the degree of~$Y_A(\u,\w)$. 
\end{proposition}
\begin{proof}
Let $r$ be the dimension of $Y_A(\u,\w)$ and its closure. The ideal $\langle \widehat{g}_{\v}(\u,\w) \rangle$ cuts out $\overline{Y_A(\u,\w)}$ plus possibly additional components $Z_1,\ldots,Z_s$ supported on $z=0$. Hence, $J$ cuts out a scheme which, set-theoretically, contains the part of $Y_A(\u,\w)$ at infinity along with possibly additional components $Z_1,\ldots,Z_s$. Since $J$ is radical, this scheme is irreducible and reduced. The dimension of $Y^{\infty}_A(\u,\w)$ is $r-1$ and by hypothesis, the dimension of each $Z_i$ is at most $r-1$. Suppose towards a contradiction that some collection, say $Z_1,\ldots,Z_t$, have dimension $r-1$. Since $\mathcal V(J)$ is irreducible, each $Z_1,\ldots,Z_t$ is set-theoretically equal to $Y^{\infty}_A(\u,\w)$. This equality holds scheme-theoretically as well: indeed, if $t>0$, then $\mathcal V(J)$ is contained in the intersection of the two components $\overline{Y_A(\u,\w)}$ and $Z_1$ of $\mathcal V(\{\widehat g_{\v}(\u,\w)_{\v \in \mathcal C_A}\})$. Hence $\mathcal V(J)$ is contained in the singular locus of $\mathcal V(\{\widehat g_{\v}(\u,\w)_{\v \in \mathcal C_A}\})$, contradicting the supposition that the scheme cut out by $J$ is reduced. In conclusion, $\overline{Y_A(\u,\w)}$ has a unique reduced component of top dimension. The degree of $Y_A(\u,\w)$ is the degree of $\overline{Y_A(\u,\w)}$ which is realized by the degree of any hyperplane section which reduces the dimension by one. Taking that hyperplane to be $\mathcal H_\infty$ finishes the proof.
\end{proof}

We aim to use Proposition \ref{prop:NaiveHomogenizationDegree} to determine the degree of $Y_A(\u,\w)$. To do so, we need to get a handle on the ideal 
\[J=\langle \widehat{g}_{\v}(\u,\w) \mid \v \in \mathcal C_A \rangle + \langle z \rangle  = \langle \widehat{g}_{\v}(\u,\w)|_{z=0} \mid \v \in \mathcal C_A \rangle.\] We observe that $\widehat{g}_{\v}(\u,\w)|_{z=0} $ is nothing but the top-degree part of $g_{\v}(\u,\w)$. We now fully describe these top-degree parts.

\begin{lemma}
\label{lem:topdegreestructureNEW}
Suppose $\v \in \mathcal C_A$ has entries in $\{-1,0,1\}$. If  $\v$ is even, then 
\[
g_{\v}(\x,\u,\w) = \pm \left(\sum\limits_{i \in \v_+}(u_i+w_i)\x^{\v-e_i} - \sum\limits_{i \in \v_-} (u_i+w_i)\x^{\v+e_i} \right) + \textit{lower order terms in }\x.
\]
If $\v$ is odd, then 
\[
g_{\v}(\x,\u,\w) = \pm 2\x^\v + \textit{lower order terms in }\x.
\]
\end{lemma}
\begin{proof}
We consider the degree $|\v|$ part of 
\[
g_{\v}(\x,\u,\w) = \x^{\v_+}(\u+\w-\x)^{\v_-} - \x^{\v_-}(\u+\w-\x)^{\v_+}.
\]
To obtain a term of degree $|\v|$ from either half of this expression, we choose $(-\x)$ as each term from $(\u+\w-\x)^{\v_{\pm}}$. For the first half, this gives $\x^{\v_+}(-\x)^{\v_-}$ and for the second half, it is $-\x^{\v_-}(-\x)^{\v_+}$ Hence, if $|\v_-|=|\v_+|$ then these terms cancel out, which occurs if and only if $\v$ is even. Otherwise, the top degree term is $2\x^{\v}$.

When $\v$ is even, we determine the degree $|\v|-1$ term of $g_\v$. Consider the term obtained from $(\u+\w-\x)^{\v_{\pm}}$  by choosing one factor to be $(u_i+w_i)$ and the rest to be $-\x$. Then we have
\begin{align*}
\x^{\v_+}(\u+\w-\x)^{\v_-} &= \sum\limits_{i \in \v_-} (-1)^{|\v_-|-1}(u_i+w_i)\x^{\v+e_i} + \textit{ lower order terms in }\x,\\
\x^{\v_-}(\u+\w-\x)^{\v_+} &= \sum\limits_{i \in \v_+} (-1)^{|\v_+|-1}(u_i+w_i)\x^{\v-e_i} + \textit{ lower order terms in }\x.
\end{align*}
Taking their difference and multiplying by $(-1)^{|\v_-|}$ gives the result for even $\v$.
\end{proof}
We write $\mydef{h_{\v}(\x,\u,\w)}$ for the top degree part $\widehat{g}_{\v}(\u,\w)|_{z=0}$ of $g_{\v}(\u,\w)$ established in Lemma \ref{lem:topdegreestructureNEW}. When $A$ is totally unimodular with even circuits, we recognize the variety cut out by $\langle h_{\v}(\u,\w) \mid \v \in \mathcal C_A \rangle$ as a \emph{reciprocal linear space}, as explained below.

Write $\mydef{\x^{-1}}=(x_1^{-1},\ldots,x_n^{-1})$ for the coordinate-wise inverse of $\x$ whenever $\x \in (\mathbb{C}^\times)^n$. Given a linear space $\mathcal L \subseteq (\mathbb{C}^\times)^n$, its \mydef{reciprocal linear space} is 
\[
\mydef{\mathcal L^{-1}} = \{\x \in (\mathbb{C}^\times)^n \mid \x^{-1} \in \mathcal L\}.
\]
For any matrix $A$, we write $\mydef{\mathcal L_A}$ for its row space. We write $\mydef{\star}$ for the \mydef{Hadamard product} of vectors, i.e. $\x \star \y = (x_1y_1,\ldots,x_ny_n)$, and $\mydef{\textrm{diag}(\x)}$ for the diagonal matrix with diagonal $\x$. Finally, for a simple statement of our next result, we write $\mydef{\s}$ as shorthand for $\u+\w$. 

\begin{lemma}
\label{lem:reciprocalRecognitionNEW}
Let $A$ be a $d \times n$ totally unimodular matrix with even circuits and let $\u,\w \in \mathbb{C}^n$ be generic. Then we have that 
\[
\mathcal V^{\times}(\{h_{\v}(\x,\u,\w)\}_{\v \in \mathcal C_A}) = \mathcal L_{A\cdot \textbf{diag}(\s^{-1})}^{-1}=(\mathcal L_{A})^{-1} \star \s.
\]
\end{lemma}
\begin{proof}
The proof is a chain of equalities:
\begin{align*}
\mathcal V^\times(\{h_{\v}(\x,\u,\w)\}_{\v \in \mathcal C_A})  &= \left\{\x \in (\mathbb{C}^\times)^n \mid \sum\limits_{i \in \v_+}s_i\x^{\v-i}-\sum\limits_{j \in \v_-}s_j\x^{\v-j} = 0 \text{ for all } \v \in \mathcal C_A\right\} \\
&= \left\{ \x \in (\mathbb{C}^\times)^n \mid \sum\limits_{i \in \v_+} s_ix_i^{-1} - \sum\limits_{j \in \v_-} s_jx_j^{-1} = 0 \text{ for all } \v \in \mathcal C_A\right\} \\
&=\left\{\x^{-1} \in (\mathbb{C}^\times)^n \mid \sum\limits_{i \in \v_+} s_ix_i - \sum\limits_{j \in \v_-} s_jx_j = 0 \text{ for all } \v \in \mathcal C_A\right\} \\
&=\left\{\x^{-1} \in (\mathbb{C}^\times)^n \mid \langle \x, \v \star \s \rangle  = 0 \text{ for all } \v \in \mathcal C_A\right\} \\
&=\left(\left(\textrm{span}\left(\left\{\v \star \s\right\}_{\v \in \mathcal C_A}\right)\right)^{\perp}\right)^{-1}\\
&= \mathcal L_{A \cdot \textbf{diag}(\s^{-1})}^{-1}.
\end{align*}
To see this last equality, we prove $\v \in \ker(A)$ if and only if $\v \star \s \in \ker(A\cdot \textbf{diag}(\s^{-1}))$. Writing columns of $A$ as $a_1,\ldots,a_n$, we have $\v \in \ker(A)$ if and only if $a_1v_1+\cdots+a_nv_n =0$. This condition is equivalent to  $(a_1s_1^{-1})(v_1s_1)+\cdots+(a_ns_n^{-1})(v_ns_n)=0$, which means $\v \star \s \in \ker(A\cdot \textbf{diag}(\s^{-1}))$.

Finally, suppose $q \in \mathcal L_{A \cdot \textbf{diag}(\s^{-1})}^{-1}$ so that $q=p^{-1}$ for $p \in \mathcal L_{A \cdot \textbf{diag}(\s^{-1})}$. Equivalently, $p = \s^{-1} \star \eta$ for $\eta \in \mathcal L_A$. This means $q = \s \star \eta^{-1}$ so $\mathcal L_{A \cdot \textbf{diag}(\s^{-1})}^{-1} \subseteq (\mathcal L_A)^{-1} \star \s$. Equality of dimension implies equality of varieties.
\end{proof}

\begin{lemma}[Proudfoot Speyer \cite{ProudfootSpeyer}]
\label{lem:proudfootspeyerNEW}
Let $A$ be a $d \times n$ totally unimodular matrix with even circuits and $\u,\w \in \mathbb{C}^n$ be generic. Then the projective variety
\[
\mathcal V(\{h_{\v}(\x,\u,\w)\}_{\v \in \mathcal C_A})
\]
is irreducible and reduced of dimension $\textrm{rank}(A)-1$ and has degree equal to $\tau_{\mathcal M(A)}(1,0)$.  The polynomials $\{h_{\v}(\x,\u,\w)\}_{\v \in \mathcal C_A}$ form a universal Gr\"obner basis for the ideal they generate.
\end{lemma}
\begin{proof}
This is a restatement of Proudfoot and Speyer's results on reciprocal linear spaces applied to our setting. We observe that the matroid underlying $A$ is the same as that underlying $A \cdot \textbf{diag}(\s^{-1})$ for a generic value of $\s \in \mathbb{C}^{n}$. 
\end{proof}

The final hypothesis of Proposition \ref{prop:NaiveHomogenizationDegree} which needs to be verified is the dimension of $Y_A(\u,\w)$. Its expected dimension is $\textrm{rank}(A)$, which is verified below to be its true generic dimension.

\begin{lemma}
\label{lem:dimYNEW}
Let $A$ be a rank $r$ totally unimodular integer matrix. For generic values of $\u$ and $\w$, every component of $X_{\Lambda(A)} \cap \{\x+\y=\u+\w\}\subseteq \mathbb{C}^{2n}$  and $Y_{A}(\u,\w)$ has dimension $r$.
\end{lemma}
\begin{proof}
Since $A$ has rank $r$ the dimension of $X_{\Lambda(A)}$ is equal to $\dim(\Lambda(A)) = r+n$. Consider the closure of $X_{\Lambda(A)} \subseteq (\mathbb{C}^\times)^{2n}$ in $\mathbb{P}^{2n}$ (that dimension is not a typo) with homogenizing variable $z$. The linear space $\{\x+\y=(\u+\w)z\}$ intersects $\overline{X_{\Lambda(A)}}$ in a union of varieties, each of projective dimension at least $r+n-n=r$. Those which are supported on the affine space $\mathbb{C}^{2n} \hookrightarrow \mathbb{P}^{2n}$ thus have affine dimension at least $r$.

On the other hand, if any component of the isomorphic variety $Y_A(\u,\w)$ has dimension greater than $r$, then the part $\overline{Y_A(\u,\w)} \cap H_{\infty}$ at infinity would have a component of projective dimension at least $r$. But $\overline{Y_A(\u,\w)} \cap H_{\infty}$ is contained in $\mathcal V(\{h_{\v}(\x,\u,\w)\}_{\v \in \mathcal C_A})$ which has affine dimension $r$ and projective dimension $r-1$ by Lemma \ref{lem:proudfootspeyerNEW}.
We deduce that any component of $Y_{A}(\u,\w) \subseteq \mathbb{C}^n$ has dimension $r$. 
\end{proof}

\begin{theorem}
\label{thm:degOfYNEW}
Suppose $A$ is totally unimodular and has only even circuits. Then
\[
\overline{Y_{A}(\u,\w)} \cap \mathcal H_{\infty} =\mathcal L_{A \cdot \textbf{diag}(\s^{-1})}^{-1}.
\]
In particular, $Y_{A}(\u,\w)$ has degree $\tau_{\mathcal M(A)}(1,0)$. 
\end{theorem}
\begin{proof}
Applying Proposition~\ref{prop:NaiveHomogenizationDegree}, we have recognized the n\"aive homogenization as the reciprocal linear space $\mathcal L_{A\cdot \textbf{diag}(\s^{-1})}$ by Lemma \ref{lem:reciprocalRecognitionNEW}. Lemma \ref{lem:proudfootspeyerNEW} gives its projective dimension, degree, and states that it is a nonreduced irreducible scheme. Lemma \ref{lem:dimYNEW} confirms the final hypothesis of Proposition~\ref{prop:NaiveHomogenizationDegree}: that the dimension of $Y_A(\u,\w)$ is one more than the dimension of its part at infinity. The theorem statement is then a restatement of the conclusion of Proposition \ref{prop:NaiveHomogenizationDegree}. 
\end{proof}

To finish the proof of our main result, we argue that the linear space $A\x = A\u$ is generic with respect to $Y_A(\u,\w)$ in the sense that they intersect in $\deg(Y_A(\u,\w))$-many points generically. Equivalently, their projective closures do not intersect at infinity. 

\begin{theorem}
\label{thm:nothingAtInfinityNEW}
Suppose $A$ is totally unimodular and has only even circuits. Then $Y_A(\u,\w)$ and $A\x=A\u$ do not intersect at infinity.
\end{theorem}
\begin{proof}
Recall that the part of $Y_A(\u,\w)$ at infinity is 
\[
\overline{Y_A(\u,\w)}\cap H_{\infty} = \mathcal L_{A \cdot \textbf{diag}(\s^{-1})}^{-1},
\]
and observe that the part of $A\x=A\u$ at infinity consists of solutions to $A\x=0$. Equivalently, this is $\mathcal L_{A}^\perp=\textrm{null}(A)$. We claim that $\mathcal L_{A \cdot \textbf{diag}(\s^{-1})}^{-1} \cap \mathcal L_{A}^\perp = \emptyset$. Suppose $q$ is in this intersection and write $q=\s\star \eta^{-1}$ for $\eta \in \mathcal L_{A}$. Note that 
\[
\langle q, \eta \rangle  = \langle s \star \eta^{-1},\eta \rangle = \sum\limits_{i=1}^n s_i\eta_i^{-1} \eta_i = \sum\limits_{i=1}^n s_i = \sum\limits_{i=1}^{n} u_i+w_i \neq 0
\]
and so $q \not\perp \eta \in \mathcal L_{A}$ and so $q \not\in \mathcal L_{A}^\perp$.
\end{proof}

\begin{theorem}[Main Theorem]
\label{thm:mlDegMainTheoremNEW}
Suppose $A$ is totally unimodular with only even circuits. Then the maximum likelihood degree of $X_{\Lambda(A)}$ is equal to $\tau_{\mathcal M(A)}(1,0)$, where $\tau_{\mathcal M(A)}$ is the Tutte polynomial of the matroid underlying $A$.
\end{theorem}
\begin{proof}
Suppose $\u$ and $\w$ are generic in $\mathbb{C}^n$. The maximum likelihood degree of $X_{\Lambda(A)}$ is the number of intersection points in $\mathbb{C}^n$, counted with multiplicity, of the variety $Y_A(\u,\w) = \mathcal V(\{g_{\v}(\x,\u,\w)\}_{\v \in \mathcal C_A})$ and the linear space $A\x = A\u$. These two varieties do not intersect at infinity, and hence, they must only intersect in points in affine space (a curve of intersections would produce a point at infinity). In particular, the linear space $A\x=A\u$ is generic in the \emph{sense of dimension} since $Y_A(\u,\w)$ has dimension $r$ and $A\x=A\u$ has codimension $r$. Even more, $A\x=A\u$ is generic in the \emph{sense of degree} since there is no intersection at infinity.  Thus the ML degree of $X_{\Lambda(A)}$ is equal to the degree of $Y_A(\u,\w)$ which we show, in Theorem \ref{thm:degOfYNEW}, to be $\tau_{\mathcal M(A)}(1,0)$.  
\end{proof}

\begin{remark}
Note that  \Cref{thm:mlDegMainTheoremNEW} holds when $A$ is totally unimodular and contains the vector of all ones in its rowspan (so, the variety $X_A$ is projective). This is because all circuits in $A$ have equal number of negative and positive entries, making the support of even size. 
Computations suggest that Theorem
\ref{thm:nothingAtInfinityNEW} always holds, even when $A$ is not totally unimodular or contains odd circuits. However, in these cases, it is unclear how to obtain the analogue of  Theorem~\ref{thm:degOfYNEW}  which relies heavily on recognizing the part of $Y_A(\u,\w)$ at infinity as a reciprocal linear space. 
\end{remark}

A main concern in applications is to classify algebraic varieties with maximum likelihood degree equal to one. When this is the case, the the maximum likelihood estimate is a rational function of the data. The varieties of the form $X_{\Lambda(A)}$ have this property only when $\mathcal M(A)$ has only circuits of size two; one may think of this condition as $\mathcal M(A)$ looking like a forest with multiple edges (see Section \ref{sec:graphs}).

\begin{lemma}
\label{lem:MobiusOne}
Let $\mathcal M$ be a matroid. Then $\tau_{\mathcal M}(1,0)=1$ if and only if $\mathcal M$ only has circuits of size $2$.
\end{lemma}
\begin{proof}
Having only circuits of size two is equal to having no loops, and no circuit of size at least $3$.

If $M$ has a loop $e$, then $\tau_M(1,0)=0$ since $\tau_M(a,b)=b\tau_{M-e}(a,b)$ which evaluates to $0$ at $b=0$. If $e$ is a coloop, then $\tau_{M}(1,0) = \tau_{M/e}(1,0) = \tau_{M-e}(1,0)$ since evaluation at $(1,0)$ counts external activity zero bases, and the complement of any such basis $B$ in $M$ is the same as the complement of~$B-e$ in $M-e$. This means that $\tau_M(1,0) \geq \tau_{M-S}(1,0)$ for any subset $S$ of $E$. Deleting everything edge not in the circuit $C$ of size $k \geq 3$ tells us $\tau_M(1,0) \geq\tau_C(1,0)$. The value $\tau_C(1,0)$ for a circuit of size $k$ is $k-1\geq 2$. 

Conversely, suppose $M$ has no loop and no circuit of size at least three. Thus there are only circuits of size two. In each circuit, there is exactly one choice for which edge may be included in a basis of external activity zero, so $\tau_M(1,0)=1$. 
\end{proof}

\begin{proposition}
\label{prop:MLdeg1}
Let $A \in \mathbb{Z}^{d \times n}$ be a totally unimodular matrix with even circuits. Then \[\textrm{mldeg}(X_{\Lambda(A)}) =1 \iff A \text{ only has circuits of size two.}\]
\end{proposition}
\begin{proof}
Since $A$ has even circuits, this condition is equivalent to $\mathcal M(A)$ having no loops and no circuits of size greater than or equal to three. Applying Lemma \ref{lem:MobiusOne} finishes the proof. 
\end{proof}

\section{Graphs}
\label{sec:graphs}

\subsection{Incidence matrices of graphs}
Fix an undirected multigraph $\mydef{G}=(\mydef{V},\mydef{E})$. The vertex-edge incidence matrix of  $G$  is the $|V|\times |E|$ integer matrix $\mydef{A_{G}}$ whose $(v,e)$-th entry is $1$ whenever~$v \in e$, and $0$ otherwise.
Associated to $G$ are two matroids, the matroid $\mathcal{M}(A_G)$ of the incidence matrix $A_G$  and the matroid $\mydef{\mathcal{M}(G)}$, called the \mydef{graphic matroid}, whose circuits are the cycles in $G$. Equivalently, the independent sets of $\mathcal M(G)$ are the forests of $G$. We remark that~$\mathcal M(A_G)$ may be distinct from $\mathcal M(G)$, as illustrated in the following example.
\begin{example}\label{ex:4graphs}
Consider the following four graphs.
\begin{figure}[!htpb]
    \centering
    \begin{tikzpicture}[
vertex/.style = {circle, fill, inner sep=2pt, outer sep=0pt},
every edge quotes/.style = {auto=left, sloped, font=\scriptsize, inner sep=1pt}
                        ]
                        
\node[] at (1,0.5) {$G_1$};
\node[vertex] (2) at (1,3)     [label=above: 2] {};
\node[vertex] (1) at (0,1)     [label=left: 1]  {};
\node[vertex] (3) at (2,1)       [label=right: 3] {};
\path   (1) edge [pos=.5, "${12}$"] (2)
        (1) edge [pos=.5, "${13}$"] (3)
        (2) edge [pos=.5, "${23}$"]  (3);
    \end{tikzpicture}  \begin{tikzpicture}[
vertex/.style = {circle, fill, inner sep=2pt, outer sep=0pt},
every edge quotes/.style = {auto=left, sloped, font=\scriptsize, inner sep=1pt}
                        ]
                        
\node[] at (1,0.5) {$G_2$};
\node[vertex] (1) at (0,3)     [label=left: 1] {};
\node[vertex] (2) at (0,1)     [label=left: 2]  {};
\node[vertex] (3) at (2,3)     [label=right: 1]  {};
\node[vertex] (4) at (2,1)       [label=right: 2] {};

\path   (1) edge [pos=.65, "${11}$"] (3)
        (1) edge [pos=.65, "${12}$"] (4)
        (2) edge [pos=.65, "${22}$"]  (4)
        (2) edge [pos=.65, "${21}$"]  (3);
    \end{tikzpicture}
\begin{tikzpicture}[
vertex/.style = {circle, fill, inner sep=2pt, outer sep=0pt},
every edge quotes/.style = {auto=left, sloped, font=\scriptsize, inner sep=1pt}
                        ]
                        
\node[] at (1,0.5) {$G_3$};
\node[vertex] (1) at (0,4)     [label=left: 1] {};
\node[vertex] (2) at (0,2)     [label=left: 2]  {};
\node[vertex] (3) at (2,5)     [label=right: 1]  {};
\node[vertex] (4) at (2,3)       [label=right: 2] {};
\node[vertex] (5) at (2,1)   [label=right: 3] {};

\path   (1) edge [pos=.65, "${11}$"] (3)
        (1) edge [pos=.65, "${12}$"]  (4)
        (1) edge [pos=.65, "${13}$"]  (5)
        (2) edge [pos=.65, "${21}$"]  (3);
    \end{tikzpicture}   \begin{tikzpicture}[
vertex/.style = {circle, fill, inner sep=2pt, outer sep=0pt},
every edge quotes/.style = {auto=left, sloped, font=\scriptsize, inner sep=1pt}
                        ]
                        
\node[] at (1,0.5) {$G_4$};
\node[vertex] (1) at (0,4)     [label=left: 1] {};
\node[vertex] (2) at (0,2)     [label=left: 2]  {};
\node[vertex] (3) at (2,5)     [label=right: 1]  {};
\node[vertex] (4) at (2,3)       [label=right: 2] {};
\node[vertex] (5) at (2,1)   [label=right: 3] {};

\path   (1) edge [pos=.65, "${11}$"] (3)
        (1) edge [pos=.65, "${12}$"]  (4)
        (1) edge [pos=.65, "${13}$"]  (5)
        (2) edge [pos=.65, "${21}$"]  (3)
        (2) edge [pos=.65, "${22}$"]  (4)
        (2) edge [pos=.65, "${23}$"]  (5);
    \end{tikzpicture}
 \caption{A three cycle, a four cycle, a tree, and the complete bipartite graph $K(2,3)$.}
    \label{fig:graphs}
\end{figure}
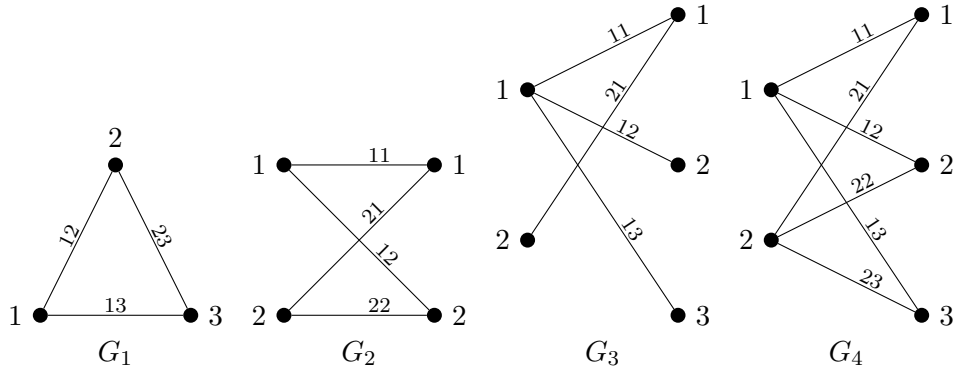

Their incidence matrices are
\begin{align*}   \scalemath{0.8}{ A_{G_1}=\begin{blockarray} {*{4}{c}} & \color{gray}{12} &  \color{gray}{13} &  \color{gray}{23} \\
\begin{block}{c[*{3}{c}]}
           \color{gray}{1} & 1 & 0 &1 \\   \color{gray}{2} &1 & 1 & 0\\   \color{gray}{3} &0 & 1 &1\\
\end{block}\end{blockarray}, \  
A_{G_2}=\begin{blockarray} {*{5}{c}} & \color{gray}{11} &  \color{gray}{12} &  \color{gray}{21} &  \color{gray}{22}  \\
\begin{block}{c[*{4}{c}]}
           \color{gray}{1} & 1 & 1 &0 & 0 \\   \color{gray}{2} &0 &0 & 1 & 1\\   \color{gray}{1} &1 &0 & 1 &0\\ \color{gray}{2} &0 &1 & 0 &1\\
\end{block}\end{blockarray}, \ A_{G_3}=\begin{blockarray} {*{5}{c}} & \color{gray}{11} &  \color{gray}{12} &  \color{gray}{13} &   \color{gray}{21}   \\
\begin{block}{c[*{4}{c}]}
           \color{gray}{1} & 1 & 1 &1 & 0 \\   \color{gray}{2} &0 &0 & 0 & 1\\ \color{gray}{1} &1 &0 & 0 &1\\ \color{gray}{2} &0 &1 & 0 &0\\  \color{gray}{3} &0 &0 & 1 &0\\  
\end{block}\end{blockarray},\ 
 A_{G_3}=\begin{blockarray} {*{7}{c}} & \color{gray}{11} &  \color{gray}{12} &  \color{gray}{13} &   \color{gray}{21}&   \color{gray}{22}  &   \color{gray}{23}  \\
\begin{block}{c[*{6}{c}]}
           \color{gray}{1} & 1 & 1 &1 & 0 &0 & 0\\   \color{gray}{2} &0 &0 & 0 & 1&1 &1\\ \color{gray}{1} &1 &0 & 0 &1 &0 &0\\ \color{gray}{2} &0 &1 & 0 &0 &1 &0\\  \color{gray}{3} &0 &0 & 1 &0 &0 &1 \\  
\end{block}\end{blockarray}}.
\end{align*}
The circuits associated to each matrix are
\begin{align*}
&\mathcal{C}_{A_{G_1}} =\emptyset, \  \mathcal{C}_{A_{G_2}} =\pm\{(1,-1,1,-1)\}, \  \mathcal{C}_{A_{G_3}} =\emptyset,\\ &\mathcal{C}_{A_{G_4}} =\pm\{(1,-1,0,1,-1,0), (1,0,-1,1,0,-1), (0,1,-1,0,-1,1)\}.
\end{align*}
We remark that the supports of these circuits correspond to the circuits of the associated graphic matroid $\mathcal M(G_i)$ for all but $i=1$. There are no circuits of $\mathcal M(A_{G_1})$, but there is a cycle in the graph $G_1$, and so $\mathcal M(G_1)$ does have a circuit.
\end{example}

A consequence of the following theorem is a characterization of those graphs $G$ for which~$\mathcal M(A_G) = \mathcal M(G)$: such $G$ are precisely the \emph{bipartite graphs}.
 \begin{theorem}\label{thm:graph_circuits}
 [\cite{Ohsugi99,villarreal1995rees},
 see also \cite[Theorem 2.3]{Petrovic14}]
Given a graph $G$, the circuits of $A_G$ are, up to multiplication by $-1$, in one-to-one correspondence with  
\begin{itemize}
\item even cycles in $G$,
\item two odd cycles in $G$ sharing a vertex,
\item two disjoint odd cycles in $G$ connected by a path.
\end{itemize}
\end{theorem}

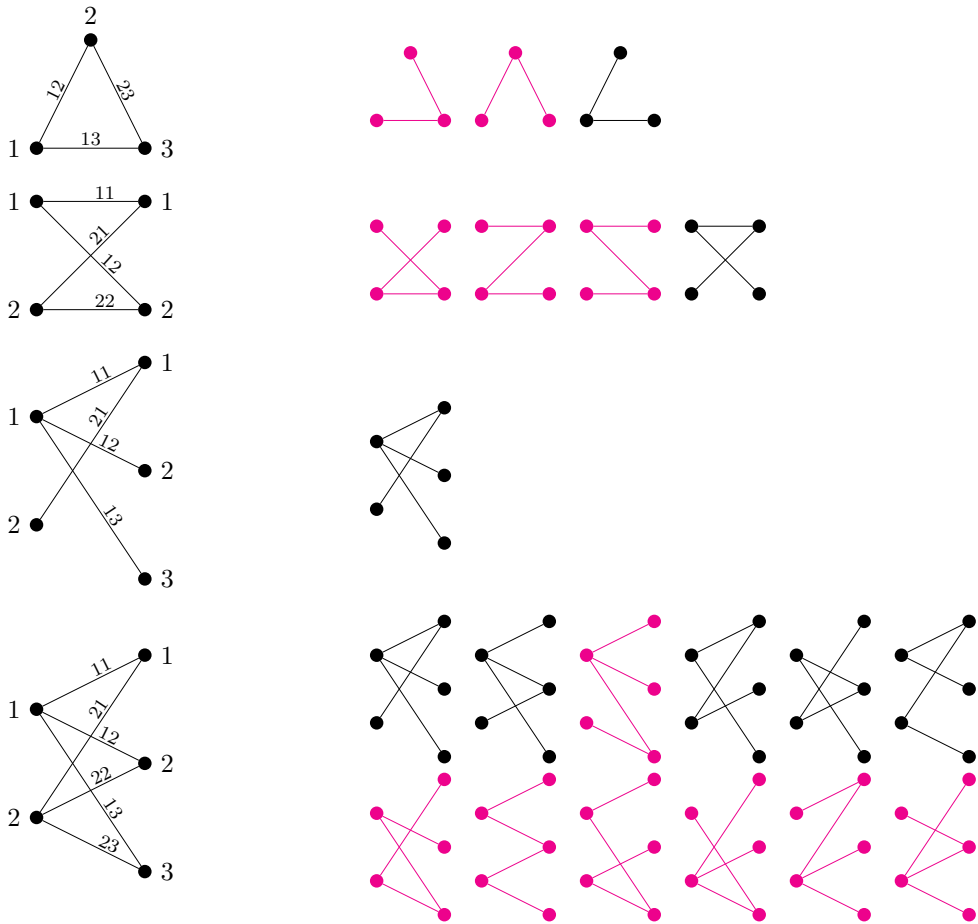
\begin{figure}[!htpb]\scalebox{0.9}{
        \centering
        \begin{tabular}{lllllll}
     
\parbox[c]{0.3\textwidth}{ 

    \begin{tikzpicture}[
vertex/.style = {circle, fill, inner sep=2pt, outer sep=0pt},
every edge quotes/.style = {auto=left, sloped, font=\scriptsize, inner sep=1pt}, scale=0.8
                        ]
                        
\node[] at (1,0.5) {};
\node[vertex] (2) at (1,3)     [label=above: 2] {};
\node[vertex] (1) at (0,1)     [label=left: 1]  {};
\node[vertex] (3) at (2,1)       [label=right: 3] {};
\path   (1) edge [pos=.5, "${12}$"] (2)
        (1) edge [pos=.5, "${13}$"] (3)
        (2) edge [pos=.5, "${23}$"]  (3);
    \end{tikzpicture}  
}          &  
            \begin{tabular}{cccccc}
 \begin{tikzpicture}[
vertex/.style = {circle, fill, inner sep=2pt, outer sep=0pt,color=magenta},
every edge quotes/.style = {auto=left, sloped, font=\scriptsize, inner sep=1pt}
,scale=0.5]
\node[vertex] (2) at (1,3)      {};
\node[vertex] (1) at (0,1)       {};
\node[vertex] (3) at (2,1)        {};

\path   (1) edge [color=magenta,pos=.5,""] (3)
        (2) edge [color=magenta,pos=.5,""]  (3);
    \end{tikzpicture}
&\begin{tikzpicture}[
vertex/.style = {circle, fill, inner sep=2pt, outer sep=0pt,color=magenta},
every edge quotes/.style = {auto=left, sloped, font=\scriptsize, inner sep=1pt}
,scale=0.5]
\node[vertex] (2) at (1,3)      {};
\node[vertex] (1) at (0,1)       {};
\node[vertex] (3) at (2,1)        {};

\path   (1) edge [,color=magenta,pos=.5,""] (2)
        (2) edge [,color=magenta,pos=.5,""]  (3);
    \end{tikzpicture}
& \begin{tikzpicture}[
vertex/.style = {circle, fill, inner sep=2pt, outer sep=0pt},
every edge quotes/.style = {auto=left, sloped, font=\scriptsize, inner sep=1pt}
,scale=0.5]
\node[vertex] (2) at (1,3)      {};
\node[vertex] (1) at (0,1)       {};
\node[vertex] (3) at (2,1)        {};

\path   (1) edge [pos=.5,""] (2)
        (1) edge [pos=.5,""] (3);
    \end{tikzpicture}
    \end{tabular}
    \\

\parbox[c]{0.3\textwidth}{ 
 \begin{tikzpicture}[
vertex/.style = {circle, fill, inner sep=2pt, outer sep=0pt},
every edge quotes/.style = {auto=left, sloped, font=\scriptsize, inner sep=1pt}, scale=0.8
                        ]
                        
\node[] at (1,0.5) {};
\node[vertex] (1) at (0,3)     [label=left: 1] {};
\node[vertex] (2) at (0,1)     [label=left: 2]  {};
\node[vertex] (3) at (2,3)     [label=right: 1]  {};
\node[vertex] (4) at (2,1)       [label=right: 2] {};

\path   (1) edge [pos=.65, "${11}$"] (3)
        (1) edge [pos=.65, "${12}$"] (4)
        (2) edge [pos=.65, "${22}$"]  (4)
        (2) edge [pos=.65, "${21}$"]  (3);
    \end{tikzpicture}
}    &
            \begin{tabular}{cccccc}
            \begin{tikzpicture}[
vertex/.style = {circle, fill, inner sep=2pt, outer sep=0pt,color=magenta},
every edge quotes/.style = {auto=left, sloped, font=\scriptsize, inner sep=1pt}
,scale=0.5]
\node[vertex] (1) at (0,3)      {};
\node[vertex] (2) at (0,1)       {};
\node[vertex] (3) at (2,3)       {};
\node[vertex] (4) at (2,1)        {};

\path   (1) edge [color=magenta,pos=.65,""] (4)
        (2) edge [color=magenta,pos=.65,""]  (4)
        (2) edge [color=magenta,pos=.65,""]  (3);
    \end{tikzpicture}
& \begin{tikzpicture}[
vertex/.style = {circle, fill, inner sep=2pt, outer sep=0pt,color=magenta},
every edge quotes/.style = {auto=left, sloped, font=\scriptsize, inner sep=1pt}
,scale=0.5]
\node[vertex] (1) at (0,3)      {};
\node[vertex] (2) at (0,1)       {};
\node[vertex] (3) at (2,3)       {};
\node[vertex] (4) at (2,1)        {};

\path   (1) edge [color=magenta,pos=.65,""] (3)
        (2) edge [color=magenta,pos=.65,""]  (4)
        (2) edge [color=magenta,pos=.65,""]  (3);
    \end{tikzpicture}
&
 \begin{tikzpicture}[
vertex/.style = {circle, fill, inner sep=2pt, outer sep=0pt,color=magenta},
every edge quotes/.style = {auto=left, sloped, font=\scriptsize, inner sep=1pt}
,scale=0.5]
\node[vertex] (1) at (0,3)      {};
\node[vertex] (2) at (0,1)       {};
\node[vertex] (3) at (2,3)       {};
\node[vertex] (4) at (2,1)        {};

\path   (1) edge [color=magenta,pos=.65,""] (3)
        (1) edge [color=magenta,pos=.65,""] (4)
        (2) edge [color=magenta,pos=.65,""]  (4);
    \end{tikzpicture}
    &
    \begin{tikzpicture}[
vertex/.style = {circle, fill, inner sep=2pt, outer sep=0pt},
every edge quotes/.style = {auto=left, sloped, font=\scriptsize, inner sep=1pt}
,scale=0.5]
\node[vertex] (1) at (0,3)      {};
\node[vertex] (2) at (0,1)       {};
\node[vertex] (3) at (2,3)       {};
\node[vertex] (4) at (2,1)        {};

\path   (1) edge [pos=.65,""] (3)
        (1) edge [pos=.65,""] (4)
        (2) edge [pos=.65,""]  (3);
    \end{tikzpicture}
    \end{tabular}
    \\

\parbox[c]{0.3\textwidth}{ 
\begin{tikzpicture}[
vertex/.style = {circle, fill, inner sep=2pt, outer sep=0pt},
every edge quotes/.style = {auto=left, sloped, font=\scriptsize, inner sep=1pt}, scale=0.8
                        ]
                        
\node[] at (1,0.5) {};
\node[vertex] (1) at (0,4)     [label=left: 1] {};
\node[vertex] (2) at (0,2)     [label=left: 2]  {};
\node[vertex] (3) at (2,5)     [label=right: 1]  {};
\node[vertex] (4) at (2,3)       [label=right: 2] {};
\node[vertex] (5) at (2,1)   [label=right: 3] {};

\path   (1) edge [pos=.65, "${11}$"] (3)
        (1) edge [pos=.65, "${12}$"]  (4)
        (1) edge [pos=.65, "${13}$"]  (5)
        (2) edge [pos=.65, "${21}$"]  (3);
    \end{tikzpicture}   

} & 
\begin{tabular}{cc}
\begin{tikzpicture}[
vertex/.style = {circle, fill, inner sep=2pt, outer sep=0pt},
every edge quotes/.style = {auto=left, sloped, font=\scriptsize, inner sep=1pt},
scale=0.5                        ]
                        
\node[vertex] (1) at (0,4)     {};
\node[vertex] (2) at (0,2)     {};
\node[vertex] (3) at (2,5)      {};
\node[vertex] (4) at (2,3)  {};
\node[vertex] (5) at (2,1)   {};

\path   (1) edge [pos=.65, ""] (3)
        (1) edge [pos=.65, ""]  (4)
        (1) edge [pos=.65, ""]  (5)
        (2) edge [pos=.65, ""]  (3);
    \end{tikzpicture}
\\
\end{tabular}    
        \\
    
\parbox[c]{0.3\textwidth}{ \begin{tikzpicture}[
vertex/.style = {circle, fill, inner sep=2pt, outer sep=0pt},
every edge quotes/.style = {auto=left, sloped, font=\scriptsize, inner sep=1pt}, scale=0.8
                        ]
                        
\node[] at (1,0.5) {};
\node[vertex] (1) at (0,4)     [label=left: 1] {};
\node[vertex] (2) at (0,2)     [label=left: 2]  {};
\node[vertex] (3) at (2,5)     [label=right: 1]  {};
\node[vertex] (4) at (2,3)       [label=right: 2] {};
\node[vertex] (5) at (2,1)   [label=right: 3] {};

\path   (1) edge [pos=.65, "${11}$"] (3)
        (1) edge [pos=.65, "${12}$"]  (4)
        (1) edge [pos=.65, "${13}$"]  (5)
        (2) edge [pos=.65, "${21}$"]  (3)
        (2) edge [pos=.65, "${22}$"]  (4)
        (2) edge [pos=.65, "${23}$"]  (5);
        \end{tikzpicture}}
            &
            
            \begin{tabular}{cccccc}
  \begin{tikzpicture}[
vertex/.style = {circle, fill, inner sep=2pt, outer sep=0pt},
every edge quotes/.style = {auto=left, sloped, font=\scriptsize, inner sep=1pt}
,scale=0.5]
\node[vertex] (1) at (0,3)      {};
\node[vertex] (2) at (0,1)       {};
\node[vertex] (3) at (2,4)       {};
\node[vertex] (4) at (2,2)        {};
\node[vertex] (5) at (2,0)    {};

\path   (1) edge [pos=.65,""] (3)
        (1) edge [pos=.65,""]  (4)
        (1) edge [pos=.65,""]  (5)
        (2) edge [pos=.65,""]  (3);
    \end{tikzpicture} & \begin{tikzpicture}[
vertex/.style = {circle, fill, inner sep=2pt, outer sep=0pt},
every edge quotes/.style = {auto=left, sloped, font=\scriptsize, inner sep=1pt}
,scale=0.5]
\node[vertex] (1) at (0,3)      {};
\node[vertex] (2) at (0,1)       {};
\node[vertex] (3) at (2,4)       {};
\node[vertex] (4) at (2,2)        {};
\node[vertex] (5) at (2,0)    {};

\path   (1) edge [pos=.65,""] (3)
        (1) edge [pos=.65,""]  (4)
        (1) edge [pos=.65,""]  (5)
        (2) edge [pos=.65,""]  (4);
    \end{tikzpicture}
&\begin{tikzpicture}[
vertex/.style = {circle, fill, inner sep=2pt, outer sep=0pt,color=magenta},
every edge quotes/.style = {auto=left, sloped, font=\scriptsize, inner sep=1pt}
,scale=0.5]
\node[vertex] (1) at (0,3)      {};
\node[vertex] (2) at (0,1)       {};
\node[vertex] (3) at (2,4)       {};
\node[vertex] (4) at (2,2)        {};
\node[vertex] (5) at (2,0)    {};

\path   (1) edge [color=magenta,pos=.65,""] (3)
        (1) edge [color=magenta,pos=.65,""]  (4)
        (1) edge [color=magenta,pos=.65,""]  (5)
        (2) edge [color=magenta,pos=.65,""]  (5);
    \end{tikzpicture}
&\begin{tikzpicture}[
vertex/.style = {circle, fill, inner sep=2pt, outer sep=0pt},
every edge quotes/.style = {auto=left, sloped, font=\scriptsize, inner sep=1pt}
,scale=0.5]
\node[vertex] (1) at (0,3)      {};
\node[vertex] (2) at (0,1)       {};
\node[vertex] (3) at (2,4)       {};
\node[vertex] (4) at (2,2)        {};
\node[vertex] (5) at (2,0)    {};

\path   (1) edge [pos=.65,""] (3)
        (1) edge [pos=.65,""]  (5)
        (2) edge [pos=.65,""]  (3)
        (2) edge [pos=.65,""]  (4);
    \end{tikzpicture}
              & 
 \begin{tikzpicture}[
vertex/.style = {circle, fill, inner sep=2pt, outer sep=0pt},
every edge quotes/.style = {auto=left, sloped, font=\scriptsize, inner sep=1pt}
,scale=0.5]
\node[vertex] (1) at (0,3)      {};
\node[vertex] (2) at (0,1)       {};
\node[vertex] (3) at (2,4)       {};
\node[vertex] (4) at (2,2)        {};
\node[vertex] (5) at (2,0)    {};

\path   (1) edge [pos=.65,""]  (4)
        (1) edge [pos=.65,""]  (5)
        (2) edge [pos=.65,""]  (3)
        (2) edge [pos=.65,""]  (4);
    \end{tikzpicture}
& \begin{tikzpicture}[
vertex/.style = {circle, fill, inner sep=2pt, outer sep=0pt},
every edge quotes/.style = {auto=left, sloped, font=\scriptsize, inner sep=1pt}
,scale=0.5]
\node[vertex] (1) at (0,3)      {};
\node[vertex] (2) at (0,1)       {};
\node[vertex] (3) at (2,4)       {};
\node[vertex] (4) at (2,2)        {};
\node[vertex] (5) at (2,0)    {};

\path   (1) edge [pos=.65,""] (3)
        (1) edge [pos=.65,""]  (4)
        (2) edge [pos=.65,""]  (3)
        (2) edge [pos=.65,""]  (5);
    \end{tikzpicture}
\\  \begin{tikzpicture}[
vertex/.style = {circle, fill, inner sep=2pt, outer sep=0pt,color=magenta},
every edge quotes/.style = {auto=left, sloped, font=\scriptsize, inner sep=1pt}
,scale=0.5]
\node[vertex] (1) at (0,3)      {};
\node[vertex] (2) at (0,1)       {};
\node[vertex] (3) at (2,4)       {};
\node[vertex] (4) at (2,2)        {};
\node[vertex] (5) at (2,0)    {};

\path   (1) edge [color=magenta,pos=.65,""]  (4)
        (1) edge [color=magenta,pos=.65,""]  (5)
        (2) edge [color=magenta,pos=.65,""]  (3)
        (2) edge [color=magenta,pos=.65,""]  (5);
    \end{tikzpicture}
& \begin{tikzpicture}[
vertex/.style = {circle, fill, inner sep=2pt, outer sep=0pt,color=magenta},
every edge quotes/.style = {auto=left, sloped, font=\scriptsize, inner sep=1pt}
,scale=0.5]
\node[vertex] (1) at (0,3)      {};
\node[vertex] (2) at (0,1)       {};
\node[vertex] (3) at (2,4)       {};
\node[vertex] (4) at (2,2)        {};
\node[vertex] (5) at (2,0)    {};

\path   (1) edge [color=magenta,pos=.65,""] (3)
        (1) edge [color=magenta,pos=.65,""]  (4)
        (2) edge [color=magenta,pos=.65,""]  (4)
        (2) edge [color=magenta,pos=.65,""]  (5);
    \end{tikzpicture}
              & 
\begin{tikzpicture}[
vertex/.style = {circle, fill, inner sep=2pt, outer sep=0pt,color=magenta},
every edge quotes/.style = {auto=left, sloped, font=\scriptsize, inner sep=1pt}
,scale=0.5]
\node[vertex] (1) at (0,3)      {};
\node[vertex] (2) at (0,1)       {};
\node[vertex] (3) at (2,4)       {};
\node[vertex] (4) at (2,2)        {};
\node[vertex] (5) at (2,0)    {};

\path   (1) edge [color=magenta,pos=.65,""] (3)
        (1) edge [color=magenta,pos=.65,""]  (5)
        (2) edge [color=magenta,pos=.65,""]  (4)
        (2) edge [color=magenta,pos=.65,""]  (5);
    \end{tikzpicture}
& \begin{tikzpicture}[
vertex/.style = {circle, fill, inner sep=2pt, outer sep=0pt,color=magenta},
every edge quotes/.style = {auto=left, sloped, font=\scriptsize, inner sep=1pt}
,scale=0.5]
\node[vertex] (1) at (0,3)      {};
\node[vertex] (2) at (0,1)       {};
\node[vertex] (3) at (2,4)       {};
\node[vertex] (4) at (2,2)        {};
\node[vertex] (5) at (2,0)    {};

\path   (1) edge [color=magenta,pos=.65,""]  (5)
        (2) edge [color=magenta,pos=.65,""]  (3)
        (2) edge [color=magenta,pos=.65,""]  (4)
        (2) edge [color=magenta,pos=.65,""]  (5);
    \end{tikzpicture}
&\begin{tikzpicture}[
vertex/.style = {circle, fill, inner sep=2pt, outer sep=0pt,color=magenta},
every edge quotes/.style = {auto=left, sloped, font=\scriptsize, inner sep=1pt}
,scale=0.5]
\node[vertex] (1) at (0,3)      {};
\node[vertex] (2) at (0,1)       {};
\node[vertex] (3) at (2,4)       {};
\node[vertex] (4) at (2,2)        {};
\node[vertex] (5) at (2,0)    {};

\path   (1) edge [color=magenta,pos=.65,""] (3)
        (2) edge [color=magenta,pos=.65,""]  (3)
        (2) edge [color=magenta,pos=.65,""]  (4)
        (2) edge [color=magenta,pos=.65,""]  (5);
    \end{tikzpicture}
& \begin{tikzpicture}[
vertex/.style = {circle, fill, inner sep=2pt, outer sep=0pt,color=magenta},
every edge quotes/.style = {auto=left, sloped, font=\scriptsize, inner sep=1pt}
,scale=0.5]
\node[vertex] (1) at (0,3)      {};
\node[vertex] (2) at (0,1)       {};
\node[vertex] (3) at (2,4)       {};
\node[vertex] (4) at (2,2)        {};
\node[vertex] (5) at (2,0)    {};

\path   (1) edge [color=magenta,pos=.65,""]  (4)
        (2) edge [color=magenta,pos=.65,""]  (3)
        (2) edge [color=magenta,pos=.65,""]  (4)
        (2) edge [color=magenta,pos=.65,""]  (5);
    \end{tikzpicture}
    \end{tabular}
        \end{tabular}}
        \caption{All spanning trees of the graphs in \Cref{fig:graphs}. The spanning trees in magenta have external activity zero. 
        }
        \label{fig:spanning_trees}
    \end{figure}

\begin{example}
\label{ex:externallyzerospanning}
Consider the four graphs in \Cref{fig:graphs}. Order the edges of each graph by the rule: edge $(i,j)>(k,l)$ if $i<k$ or $i=k$ and $j<k$.  \Cref{fig:spanning_trees} shows the spanning trees of each graph, with those of external activity zero coloured magenta. 
\end{example}

We now write our main result in the context of graphs.

\begin{theorem}\label{thm:mld_graphs}
Fix a graph $G$. Then $\mathcal M(A_G) = \mathcal M(G)$ if and only if $G$ is bipartite. In this case, 
 \begin{align*}
 \textrm{deg}(X_{\Lambda(A_G)}) &= \tau_G(1,1)=\text{number of spanning forests in } G,
  \\
  \textrm{mldeg}(X_{\Lambda(A_G)}) &= \tau_G(1,0)=\text{number of spanning forests with external activity zero  in }G.
  \end{align*}
In particular, $\mathrm{mldeg}(X_{\Lambda(A_G)})=1$ if and only if $G$ only has cycles of size two.
\end{theorem}

\begin{proof}
A graph $G$ is bipartite if and only if all of its cycles are even. Hence, the circuits of $G$ are in one-to-one correspondence with its cycles by \Cref{thm:graph_circuits}. That $A_G$ is totally unimodular may be proven by a simple induction argument. Evaluating the Tutte polynomial  at $(1,1)$ and $(0,1)$ yields
\begin{align*}
    &\tau_G(1,1)=\sum\limits_{i,j}t_{ij}1^i1^j=\text{number of spanning forests in } G, \text{ and }\\
    &\tau_G(1,0)=\sum\limits_{i,j}t_{ij}1^i0^j=\sum\limits_{i}t_{i0}1^i0^0=\text{number of spanning forests with external activity zero in } G.
\end{align*}
Applying \Cref{prop:MLdeg1} gives the characterization of when the ML degree is one.
\end{proof}

\begin{example}\label{ex:4graphs2}
Of the four graphs in \Cref{fig:graphs}, only $G_1$ is not bipartite. Hence, we may use \Cref{thm:mld_graphs} to obtain degrees and maximum likelihood degrees of $X_{\Lambda(A_{G_2})} ,X_{\Lambda(A_{G_3})} $, and $X_{\Lambda(A_{G_4})} $, that is, these ML degrees can be read-off from \Cref{fig:spanning_trees}.
     \begin{enumerate}
    \item $\mathrm{deg}(X_{\Lambda(A_{G_2})})=4$ and $\mathrm{mldeg}(X_{\Lambda(A_{G_2})})=3$. \\ In general, for a cycle $C_n$ on $n$ vertices, $\mathrm{deg}(X_{\Lambda(A_{C_n})})=n$ and $\mathrm{mldeg}(X_{\Lambda(A_{C_n})})=n-1$. 
    \item  $ \deg(X_{\Lambda(A_{G_3})})=1$ and $\mathrm{mldeg}(X_{\Lambda(A_{G_3})})=1$. \\ See \Cref{thm:mld_graphs}.
    \item  $\deg(X_{\Lambda(A_{G_4})})=12$ and $\mathrm{mldeg}(X_{\Lambda(A_{G_4})})=7$.\\ When $G$ is complete bipartite, see  \Cref{tab:CompleteBipartiteDegrees} and \Cref{cor:bipartite}. 
\end{enumerate}
For the non-bipartite graph $G_1$, the matrix $A_{G_1}$ is square of full-rank and so the toric variety~$X_{\Lambda(A_{G_1})}$ is the entire affine space. Thus, it has degree and maximum likelihood degree one.
\end{example}
\begin{table}[!htpb]\scalebox{0.95}{
\begin{tabular}{|c||r|r|r|r|r|r|}
\hline
&1 &2 & 3 & 4 & 5 & 6 \\ \hline \hline       
1 & 1&1&1&1&1&1\\
2&      1&4&12&32&80&192\\
3&      1&12&81&432&2\,025&8\,748\\
4&      1&32&432&4\,096&32\,000&221\,184\\
5&      1&80&2\,025&32\,000&390\,625&4\,050\,000\\
6&      1&192&8\,748&221\,184&4\,050\,000&60\,466\,176\\
 \hline
\end{tabular} \quad 
\begin{tabular}{|c||r|r|r|r|r|r|}
\hline
&1 &2 & 3 & 4 & 5 & 6 \\ \hline \hline 
1&      1&1&1&1&1&1\\
2&      1&3&7&15&31&63\\
3&      1&7&31&115&391&1\,267\\
4&      1&15&115&675&3\,451&16\,275\\
5&      1&31&391&3\,451&25\,231&164\,731\\
6&      1&63&1\,267&16\,275&164\,731&1\,441\,923 \\ \hline
\end{tabular}}
\caption{Degree (Left) and ML degree (Right) of $X_{\Lambda(A_{K(m_1,m_2)})}$ for~$m_1,m_2=1\dots,6$.}
\label{tab:CompleteBipartiteDegrees}
\end{table} 
It is known that the number of spanning trees in the complete bipartite graph is $K(m_1,m_2)$ is~$\tau_{K(m_1,m_2)}(1,1)=m_1^{m_2-1}m_2^{m_1-1}$, given as OEIS A072590. The values of $\tau_{K(m_1,m_2)}(1,0)$ are OEIS A136126. This later number is, equivalently, the number of acyclic orientations of $K(m_1,m_2)$ with a unique source \cite[Theorem 4.1]{Parking}.
We write the first of these values in  \Cref{tab:CompleteBipartiteDegrees}.

\begin{corollary}\label{cor:bipartite}
    Let $K(m_1,m_2)$ be the complete bipartite graph. Then, 
    \begin{align*}
        \deg(X_{\Lambda(A_{K(m_1,m_2)})})= m_1^{m_2-1}m_2^{m_1-1}  \text{ and } \mathrm{mldeg}(X_{\Lambda(A_{K(m_1,m_2)})})= \sum\limits_{k=1}^{m_1} \left(\frac{1}{k} \sum\limits_{i=1}^{k} (-1)^{m_1-i}{k\choose i}i^{m_1}\right)k^{m_2}.
    \end{align*}
\end{corollary}

For ease-of-reading, we specify the ML degree formula when $m_1=m$ and $m_2=2,3,4$ and $5$. That is, the following are expressions for the rows of the right table of \Cref{tab:CompleteBipartiteDegrees}.
\begin{align*}
 & \mathrm{mldeg}(X_{\Lambda(A_{K(m,2)})})=  2^m -1,\\ 
  & \mathrm{mldeg}(X_{\Lambda(A_{K(m,3)})})= 2\cdot 3^m+3\cdot 2^m +1,\\ 
  & \mathrm{mldeg}(X_{\Lambda(A_{K(m,4)})})= 6\cdot 4^m - 12\cdot 3^m + 7\cdot 2^m - 1,\\
   & \mathrm{mldeg}(X_{\Lambda(A_{K(m,5)})})= 24\cdot 5^m - 60\cdot 4^m + 50\cdot 3^m - 15\cdot 2^m + 1. 
\end{align*}
\Cref{cor:bipartite} finds direct applications in algebraic statistics, including the confirmation of conjectures in \cite{coons2020generalized} and  \cite{Amendola17} on the degree and maximum likelihood degree, respectively, of no-three-way interaction models. We explain this in Section \ref{sec:statistics}.  

\subsection{Incidence matrices of directed graphs} We now specify our main theorem to a more general setting than graphs: directed graphs. 
Let $\mydef{\overrightarrow{G}}=\mydef{(V,\overrightarrow{E})}$ be a \emph{directed graph}. The vertex-edge incidence matrix of  $\overrightarrow{G}$   is the $|V|\times |\overrightarrow{E}|$ integer matrix $A_{\overrightarrow{G}}$  whose $(u,\overrightarrow{e})$-th entry is $1$ and $(v,\overrightarrow{e})$-th entry is $-1$, whenever $\overrightarrow{e}=(u,v)$, and zero otherwise. 
\begin{figure}[h]
\includegraphics[scale=1]{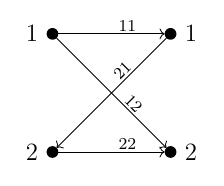} \quad \includegraphics[scale=1]{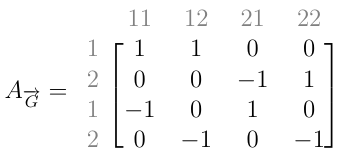} 
\caption{A directed bipartite graph and its incidence matrix. Up to sign, the matrix has only one circuit $\mathcal{C}_{A_{\vec{G}}}=\{(1,-1,1,1)\}$.}
\label{fig:directedgraph}
\end{figure}
An example of $\overrightarrow{G}$ and $A_{\overrightarrow{G}}$ is given in Figure \ref{fig:directedgraph}.
The matroid of $A_{\overrightarrow{G}}$ is widely studied and intimately connected to the notion of an \emph{oriented matroid} (see \cite{oxley2006matroid}). Any matrix of the form $A_{\overrightarrow{G}}$ is totally unimodular  (see \Cref{rmk:graphs}). Unlike the undirected case (see \Cref{thm:graph_circuits}), the circuits of $\mathcal{M}(A_{\overrightarrow{G}})$ always coincide with the cycles of $G$. Hence, when $\overrightarrow{G}$ is bipartite, its cycles are all even and so one may apply \Cref{thm:mlDegMainTheoremNEW} to obtain the following corollary.  

\begin{corollary}\label{cor:mld_oriented_graphs}
Let $\overrightarrow{G}$ be a  directed graph with underlying graph $G$.  Then 
\[
\textrm{deg}(X_{\Lambda(A_{\overrightarrow{G}})}) = \tau_G(1,1)=\text{number of spanning forests in }G.
\]
If $G$ is bipartite, then
\[
\textrm{mldeg}(X_{\Lambda(A_{\overrightarrow{G}})}) = \tau_G(1,0)=\text{number of spanning forests with external activity zero  in }G
\]
and the ML degree is one when $G$ only has cycles of size two.
\end{corollary}

\subsection{Incidence matrices of signed graphs} We consider one last graph-theoretic construction and its relationship to our results. Specifically, we address matrices which come from \emph{signed graphs}. This covers all matrices with entries $0, \pm 1$ whose columns have exactly two nonzero elements.

A \mydef{signed graph} is a triple $\mydef{G^{\pm}}=(V,E,\delta)$ where $V$ and $E$ are vertices and edges, as usual, and~$\mydef{\delta}:E \rightarrow \{\pm\}$ is a \mydef{signed function}. The vertex-edge {incidence matrix} for  $G^{\pm}$  is the $|V|\times |E|$ integer matrix $A_{G^{\pm}}$ such that for each edge $e=(u,v)\in E$, the $(u,e)$-th entry is $1$ or $-1$, subject to the condition that this entry agrees with the $(v,e)$-th entry when $\textrm{sign}(e)=-$ and is the opposite sign when $\textrm{sign}(e)=+$. We remark that this is well-defined up to the signs of each column, which makes the matroid underlying this matrix well-defined. Such a matroid is called a \mydef{frame matroid} or \mydef{signed matroid} (see, for instance, \cite{ZASLAVSKY198247}).
A signed graph where  the product over the signs in any cycle is positive is called a \mydef{balanced signed graph}.

The notion of a signed graph, as well as its associated incidence matrix, interpolates between the directed and undirected cases discussed previously: incidence matrices of directed graphs are incidence matrices of all-positive signed graphs, whereas, incidence matrices of undirected graphs correspond to the all-negative sign function.

Signed graphs with the same matroid will be useful to us. For this we need the definition of switching.
Given a signed graph $G^{\pm}=(G,\delta)$ and a function $\mydef{\varsigma}: V(G)\rightarrow \{\pm\}$, the operation of \mydef{switching} $G^{\pm}$ by $\varsigma$ results in a new signed graph $(G,\delta^{\varsigma})$ called the \mydef{switched graph}. The switched graph has the same underlying graph as $G$ but its sign function is defined by 
\begin{align*}
\mydef{\delta^{\varsigma}}:E &\to \{\pm\} \\
(u,v)=e &\mapsto \varsigma(u)\delta(e)\varsigma(v).
\end{align*} Two signed graphs are \mydef{switching equivalent} if one is obtained from another via a sequence of switches.

\begin{figure}[h]
\includegraphics[scale=1]{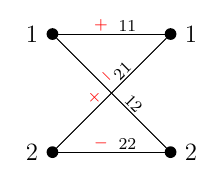} \quad \includegraphics[scale=1]{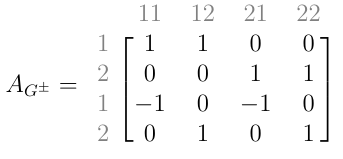}
\caption{A directed bipartite graph and its incidence matrix. Up to sign, the matrix has one circuit $\mathcal{C}_{A_{G^{\pm}}}=\pm \{(1,-1,-1,1)\}$.}
\end{figure}

We collect several relevant results about signed graphs in the following theorem.
\begin{theorem}\label{thm:Z_graphs}
 \cite[Corollary~3.3+Proposition~8A.5]{ZASLAVSKY198247} and   \cite[Corollary~5.4]{ZASLAVSKY198247}
 
 Let $G^\pm$ be a signed graph. The following are equivalent:
\begin{itemize}
\item $G^\pm$ is balanced,
\item $G^\pm$ is switching equivalent to an all-positive signed graph,
\item the incidence matrix of $G^\pm$ is totally unimodular.
\end{itemize} Additionally, two signed graphs with the same underlying undirected graph are switching equivalent if and only if their incidence matrices have the same frame  matroid.
\end{theorem}

\begin{remark}\label{rmk:graphs}
Since incidence matrices of signed graphs interpolate between those of directed and undirected graphs, one can directly see from \Cref{thm:Z_graphs} that the incidence matrix of a directed graph is totally unimodular, whereas the incidence matrix of an undirected graph is totally unimodular if and only if the graph is bipartite.
\end{remark}

The following corollary puts our main theorem in the context of signed graphs.

\begin{corollary}
Let $G^{\pm}$ be a balanced signed graph on underlying graph $G$.   Then
\[
\textrm{deg}(X_{A_{G^{\pm}}}) = \tau_G(1,1)=\text{number of spanning forests in } G.
\]
If $G$ is bipartite, then
\[
  \textrm{mldeg}(X_{A_{G^{\pm}}}) = \tau_G(1,0)=\text{number of spanning forests with external activity zero  in }G
  \]
  and the ML degree is one if and only if $G$ has only cycles of size two.
\end{corollary}

\section{Algebraic Statistics}
\label{sec:statistics}

Many discrete statistical models are of form $X\cap \Delta_{n-1}$ where $X$ is some algebraic variety in $\mathbb C^n$ and \mydef{$\Delta_{n-1}$} is the interior of the standard $(n-1)$-simplex in $\mathbb{R}^n$. \emph{Log-linear models} are precisely the statistical models appearing this way when $X$ is a toric variety. 
More precisely, given $A \in \mathbb{Z}^{d \times n}$, the \mydef{log-linear model} (or \mydef{toric model}) for $A$ is the set
\begin{align*}
\{\p\in \Delta_{n-1} \mid \log \p\in \mathcal L_A\}=X_A\cap \Delta_{n-1},
\end{align*} where $\mydef{\log}:\Delta_{n-1} \to \mathbb{C}^n$ is the coordinate-wise logarithm and $\mathcal L_A$ is, as before, the row space of~$A$. The columns of~$A$ are indexed by \emph{states} of the model and the rows are indexed by its \emph{parameters}.
This setup is illustrated on a practical example in \Cref{ex:Reye's Syndrome} where the toric model is of the form~$X_{\Lambda(K(2,3))} \cap \Delta_{11}$. Several other examples may be found in the literature \cite{diaconis1998algebraic,agresti2012categorical}, in particular, \cite{Sullivant.2018} uses a model of this form as a running example (Plum root cuttings) throughout the book. 

Log-linear models are widely used in applications; they include graphical models, hierarchical models, quasi-independence models, decomposable Bayesian networks, balanced staged tree models, and several group-based phylogenetic models. Many of these admit the structure of a Lawrence lift, $X_{\Lambda(A)} \cap \Delta_{2n-1}$. Among these, is the \emph{no-three-way interaction model}, which was the initial motivation for this paper. 

Given some empirically observed data and a discrete model, the problem of maximum likelihood estimation asks for the distribution on that model which best explains the data. Specifically, suppose  $\mydef{\mathbf{u}}\in \mathbb{Z}_{\geq 0}^n$ is the vector of counts for the given data with sample size $\mydef{N}=u_1+\cdots+u_{n}$. 
 The  probability of observing the data~$\mathbf{u}$, given that the underlying probability distribution $\mathbf{p}$ belongs to $X\cap \Delta_{n-1}$ (see \cite[Section~1.1]{Sullivant.2018})
is represented by the \mydef{likelihood function} 
\begin{align}\label{eq:discretelikeq}
\p \mapsto {\dfrac{N!}{\prod\limits_{i=1}^{n}(u_{i}!)}}p_{1}^{u_{1}}\cdots p_{n}^{u_{n}}.
\end{align}
The \mydef{maximum likelihood estimate} (MLE) for the data $\mathbf{u}$ in the statistical model $X\cap \Delta_{n-1}$ is the value $\hat{\p} \in X \cap \Delta_{n-1}$ which maximizes the likelihood function.
Removing the constant multinomial coefficient in \Cref{eq:discretelikeq} does not change the MLE. Similarly, since $\log$ is monotonic, applying the natural logarithm to the likelihood function does not change the MLE. Performing both of these operations, we obtain the \mydef{log-likelihood function}
$$\mydef{\ell(\mathbf{p}|\mathbf{u})}=u_{1}\mathrm{log}(p_{1})+\cdots + u_{n}\mathrm{log}(p_{n}).$$

The MLE of a model is found among the critical points of  $\ell(\mathbf{p}|\mathbf{u})$.  The maximum likelihood degree of a statistical model $X \cap \Delta_{n-1}$ is the number of complex critical points of $\ell(\mathbf{p}|\mathbf{u})$ in $X \cap \Delta_{n-1}$, for generic data $\u$. This number is well-defined and agrees with  \Cref{def:mlequationsNEW} from Section \ref{sec:mldegree}.

\begin{example}
\label{ex:Reye's Syndrome}
Reye's syndrome is a condition that causes swelling in the liver and brain. It can affect children and teenagers after they have an infection. It is suspected that using Aspirin to treat an infection increases the chances of getting Reye's syndrome \cite{belay1999reye}. One designs an  experiment to investigate the effect of infection, age, and Aspirin in patients with Reye's syndrome. The hypothetical data in \Cref{tb:data}  cross-classified according to the features of interest.
\begin{table}[!htpb]
\label{tb:data}
    \centering
    \begin{tabular}{ccccc} 
         &  &  &  \textbf{Infection} & \\ 
        \textbf{Take Aspirin Reg.} & \textbf{Age} &  Varicella&  Influenza& Gastroenteritis\\ 
        \hline
        Yes & Child & 19 & 11 & 2\\ 
         & Teen & 9  & 10 & 1\\ 
        No & Child & 400 & 100 & 74\\ 
         & Teen & 304 & 88 & 50\\ 
    \end{tabular}
    \caption{Hypothetical data $\u=(u_{111},u_{112},\cdots, u_{232})=(19,9,\cdots 50)$ on $1070$ patients with Reye's Syndrome, cross-classified according the previous infection, whether the patient used Aspirin regularly, and age. \textit{Note:} Although hypothetical, this data is inspired by real-world experiments \cite{belay1999reye}. }
    \label{tab:reye}
\end{table}

Let $X_1$ denote the binary random variable for the Aspirin intake, let $X_2$ be the random variable for type of infection, and $X_3$ be the random variable for the age group of the patient.  The state space of the experiment is  $\Omega=\Omega_1\times \Omega_2\times \Omega_3$, where
 \[\Omega_1=\left\{\begin{tabular}{c}\text{yes Aspirin}\\ \text{no Aspirin}\end{tabular}\right\},\quad  \Omega_2= \left\{\begin{tabular}{c}\text{Varicella}\\ \text{Influenza}\\ \text{Gastroenteritis}\end{tabular}\right\},   \quad \Omega_3= \left\{\begin{tabular}{c}\text{Child}\\\text{Teen}\end{tabular}\right\}.\]
Denote $\mydef{p_{ijk}}=P(X_1=i,X_2=j,X_3=k) \text{ for } (i,j,k)\in \Omega.$
A probability distribution for data~$\u$ is of form $\p=(p_{ijk})_{(i,j,k)\in \Omega}\in \Delta_{11}$. The \emph{no-three-way interaction model} for our three random variables is the set of all such probability distributions which are parametrized in such a way that all two-way relation dependencies are possible, but no relation dependency among all three exists, that is, 
\[ p_{ijk}=\dfrac{1}{N(\textbf{a},\textbf{b},\textbf{c})}a_{i,j}b_{j,k}c_{i,k} \text{ for any } (i,j,k)\in\Omega,\]
where $\textbf{a},\textbf{b},\textbf{c}$ are some potential functions and \(
 \mydef{N(\textbf{a},\textbf{b},\textbf{c})}= \sum\limits_{(i,j,k)\in \Omega}a_{i,j}b_{j,k}c_{i,k}.\)

The parametrization of the model naturally induces the familiar monomial map
\begin{align*}
{\varphi}_{\Lambda(A_{K(2,3)})}: \mathbb C[p_{ijk} \mid (i,j,k)\in [2]\times [3]\times [2]] &\rightarrow \mathbb C[a_{ij},b_{jk},c_{ik}\mid (i,j,k)\in [2]\times [3]\times [2]],\\
p_{ijk}&\mapsto a_{ij} b_{jk} c_{ik}.
\end{align*}
This monomial map comes from the matrix $\Lambda(A_{K(2,3)})$, the Lawrence lift of the vertex-edge incidence matrix of the (undirected) complete bipartite graph $K(2,3)$. 
So, this no-three-way interaction model is the log-linear model $X_{\Lambda(A_{K(2,3)})}\cap \Delta_{11}$ obtained as a Lawrence lift of $X_{A_{K(2,3)}}\cap \Delta_5$. 
From \Cref{ex:4graphs2},  the degree of the model is $12$ and its maximum likelihood degree is $7$.
\end{example}

\subsection{No three way interaction models/three dimensional transportation problem}

Following \Cref{ex:Reye's Syndrome}, a no-three-way interaction model pertains to three discrete random variables $\mydef{X_1,X_2,X_3}$ with $\mydef{m_1, m_2, m_3}$ number of states, respectively. It models the situation when all two-way relation dependencies are assumed, but  no relation dependency among  the three variables simultaneously exists. Precisely, a \mydef{no-three-way interaction model} is a model of the form
\begin{align*}
\mydef{M_{m_1,m_2,m_3}}=\left\{(p_{ijk})\in \Delta_{m_1m_2m_3-1} \mid \dfrac{1}{N(\textbf{a},\textbf{b},\textbf{c})}a_{i,j}b_{j,k}c_{i,k} \text{ for all } (i,j,k)\in [m_1]\times[m_2]\times [m_3]\right\}.
\end{align*}

Often, no-three-way interaction models appear in the literature of optimal transportation problems under the name {\emph{three dimensional transportation model}} \cite{ishizeki2003standard,sturmfels1997variation,Sturmfels96}. 
A typical application is to optimize the allocation of $m_1$ different products from $m_2$ different factories to $m_3$ different consumer zones, so as to minimize the total transportation cost; this corresponds to optimization over $M_{m_1m_2m_3}$ (see also  \cite{vlach1986conditions}).

The model $M_{m_1m_2m_3}$ is the intersection of the vanishing of the kernel of
 \begin{align*}
{\varphi}: \mathbb C[p_{ijk} \mid (i,j,k)\in [m_1] \times [m_2]\times [m_3]] & \rightarrow \mathbb C[a_{ij},b_{jk},c_{ik}\mid (i,j,k)\in [m_1]\times [m_2]\times [m_3]],\\
p_{ijk} & \mapsto a_{ij} b_{jk} c_{ik},
\end{align*}
with the probability simplex.
When one of the random variables is a binary random variable, then the matrix encoding $\varphi$ is a Lawrence lift of a complete bipartite graph on the sizes of the state spaces of the remaining two random variables. We can now apply \Cref{cor:bipartite} to obtain some of their degrees.

\begin{corollary}({Positive resolution of \cite[Conjecture~36]{Amendola17}})
\label{cor:no-two-wayinteraction}
The degrees and maximum likelihood degrees of the no-three-way interaction models (c.f. three dimensional transportation models) of the form $M_{m_1,m_2,2},M_{m_1,2,m_2},M_{2,m_1,m_2}$ are given by the formulae in \Cref{cor:bipartite}. In particular, the ML degree of $M_{m,2,2}$ is $2^m-1$, as conjectured in \cite{Amendola17}.
\end{corollary}

Despite the simple description of $M_{m_1,m_2,m_3}$, as $m_1,m_2,$ and $m_3$ increase, the ideal which cuts out $M_{m_1,m_2,m_3}$ gets surprisingly complicated (see, for instance \cite{almendra2023markov,de2006markov}).  No-three-way interaction models are discussed in several papers 
\cite{almendra2023markov, agresti2012categorical,Amendola17, bayer2001sysygies,Bernstein15, Bernstein17, coons2020generalized, de2006markov, diaconis1998algebraic} 
 and in books \cite{Sturmfels96,Sullivant.2018}. They are the smallest irreducible hierarchical models (see \Cref{ex:no3wayIsHier}) and knowledge about their ML degrees provides knowledge about  models built from them via toric fiber products (the ML degrees involved in a toric fiber product multiply~\cite{amendola2020maximum}).

\begin{remark}
For no-three-way interaction models involving no binary random variable, the matrix encoding the monomial map may be interpretted as a \emph{higher Lawrence lift} \cite{santos2003higher}. For example, the matrix associated to $M_{333}$  is the third Lawrence lift:
\[ \Lambda^{(3)}(A_{K(3,3)})=
\begin{bmatrix}
A_{K(3,3)} & 0 & 0 \\
0 & A_{K(3,3)} & 0 \\
0 & 0 & A_{K(3,3)} \\
I_9 & I_9 & I_9
\end{bmatrix}
\]
The degree of $M_{333}$ is $2457$ but have no conjecture for a combinatorial interpretation of this number. We leave the case when there is no binary random variable to future work.
\end{remark}

\subsection{Three dimensional quasi-independence models}\label{sec:n2w} 
A $k$-dimensional \mydef{independence model} for random variables $X_1,\ldots, X_k$ with $m_1,\ldots,m_k$ number of states, respectively, is a log-linear model encoded by an integer matrix with columns of the form 
\[e_{i_1} + e_{m_1+i_2}+\ldots + e_{m_1\cdots m_{k-1}+i_k}\quad \text{ for }\quad (i_1,\ldots,i_k)\in [m_1]\times \cdots \times [m_k].\] 
\mydef{Quasi-independence models} are log-linear models coming from matrices which are obtained by removing columns from a matrix giving an independence model. Statistically, this means that the removed columns correspond to states with probability zero, that is,  some combinations of states  cannot occur together, and otherwise the random variables are independent of eachother. A standard reference for quasi-independence  is \cite{aoki2012markov}. 
\begin{example}
Matrices $A_{G_2}$ and $A_{G_4}$  in \Cref{ex:4graphs} are examples of two-dimensional independence models in $(2,2)$ and $(2,3)$ number of states, respectively. The matrix $A_{G_3}$ is a quasi-independence model obtained by removing columns labeled $13, 22, 23$ from $A_{G_4}$. 
\end{example}

We apply the results in this paper to three-dimensional quasi-independence models that arrive as Lawrence lifts of two dimensional quasi-independence models. The two-dimensional quasi-independence models are determined by  incidence matrices of a bipartite graphs.

The construction of the model 
$X_{\Lambda(A_G)}\cap \Delta_{2n-1}$  and \Cref{thm:mld_graphs} yields the following application. 
\begin{corollary}[Implication of \Cref{thm:mld_graphs} for quasi-independence models]
\label{cor:applytoquasi}
Consider the three-dimensional quasi-independence model  $X_{\Lambda(A_G)}\cap \Delta_{2n-1}$ where $G$ is a simple bipartite graph.  Then 
\begin{align*}
\textrm{deg}(X_{\Lambda(A_G)}\cap \Delta_{2n-1}) &= \tau_G(1,1)=\text{number of spanning forests in } G,\\
\textrm{mldeg}(X_{\Lambda(A_G)}\cap \Delta_{2n-1}) &= \tau_G(1,0)=\text{number of spanning forests with external activity zero  in }G.
\end{align*}
  In particular,  $X_{\Lambda(A_G)}\cap \Delta_{2n-1}$  has maximum likelihood degree one if and only if~$G$ is a forest.
\end{corollary}

 Work on the  MLE of two-dimensional quasi-independence models is contained in \cite{bishop2007discrete} and more recently in  \cite{coons2021quasi}, with an emphasis on  those with  maximum likelihood degree  equal to one. 
\Cref{cor:applytoquasi} emphasizes  that the Lawrence lift complicates the model; while two dimensional quasi-independent models have rational MLE (maximum likelihood degree equal to one) whenever the graph is doubly chordal \cite{coons2021quasi}, their Lawrence lifts have rational MLE if and only if the graph is a forest.

\subsection{Hierarchical models}
\label{sec:hierarchical}
 Hierarchical models are log-linear models  that describe the dependency relations among discrete random variables. They have a monomial parametrization determined by a simplicial complex \mydef{$\Gamma$} on ground set $\mydef{[n]}=\{1,2,\ldots,n\}$ and a vector~$\mydef{\mathbf{r}}\in \mathbb{N}^n$ of positive integers. Vertex $i$ in $\Gamma$ represents a  random variable $X_i$  with $r_i$ number of states. A set of vertices is not in the same face  of $\Gamma$ if and only if the corresponding random variables are independent given that the rest of the variables are fixed.  

Before we formally define the matrix that parametrizes a hierarchical model,  we need to simplify some notation. For any $F\subseteq [n]$ let  $\mydef{r_F}= \prod\limits_{i\in F}r_i$ with the convention $r_\emptyset=1$, and $\mydef{[\r]_F}=\prod\limits_{i\in F}[r_i]$. When $F=[n]$  we will use the shorthand $\mydef{r}=r_{[n]}$ and $\mydef{[\r]}=[\r_{[n]}]$. Let $\mydef{\textrm{facet}(\Gamma)}$ be the set of facets of~$\Gamma$, that is, maximal simplices in $\Gamma$. Denote the sum of $\{r_F\}_{F \in \textrm{facet}(\Gamma)}$ by the letter $\mydef{d}$.  The notation~$\mydef{\i_F}$ is for the vector obtained by selecting only the terms indexed by $F$ in the vector $\mathbf{i}\in [\r]$.  

A hierarchical model for $\Gamma$ and $\r$ is represented by the  matrix  $A= \mydef{A_{\Gamma,\mathbf{r}}}\in \mathbb{Z}^{d\times r}$   as follows: for~$F\in \mathrm{facet}(\Gamma)$, $\mathbf{i}\in [\r]_F$,  and $\mathbf{k}\in [\r]$, define
\begin{align*}
A_{(F,\mathbf{i}),\mathbf{k}} =\begin{cases} 1  & \text{if} \ \mathbf{k}|_{F}= \mathbf{i},\\
0  & \text{else}.
\end{cases}
\end{align*}
The hierarchical model for  $\Gamma$ and $\r$  is the log-linear model $X_{A_{(\Gamma,\mathbf{r})}}\cap \Delta_{r-1}$.

\begin{example}\label{ex:marginalmatrix}
Let $\mathbf r=(3,2,2)$ and $\Gamma$ have facets $F=\{1,2\}$ and $T=\{2,3\}$. Then, the  matrix for the hierarchical model on $\Gamma,\mathbf r$ is

\begin{minipage}{6cm}
\begin{align*} & A_{\Gamma,\mathbf r}=\begin{blockarray} {*{13}{c}} & \color{gray}{111} &  \color{gray}{112} &  \color{gray}{121}  &  \color{gray}{122}& \color{gray}{211} &  \color{gray}{212} &  \color{gray}{221}  &  \color{gray}{222}& \color{gray}{311} &  \color{gray}{312} &  \color{gray}{321}  &  \color{gray}{322}\\
\begin{block}{c[*{12}{c}]}
           \color{gray}{F11} & 1 & 1 &0 & 0 &0 & 0 &0 & 0 &0 & 0 &0 & 0 \\   \color{gray}{F12} & 0 & 0 &1 & 1 &0 & 0 &0 & 0 &0 & 0 &0 & 0 \\   \color{gray}{F21} & 0 & 0 &0 & 0 &1 & 1 &0 & 0&0 & 0 &0 & 0  \\ \color{gray}{F22} & 0 & 0 &0 & 0 &0 & 0 &1 & 1 &0 & 0 &0 & 0 \\\color{gray}{F31} & 0 & 0 &0 & 0 &0 & 0 &0 & 0&1 & 1 &0 & 0  \\ \color{gray}{F32} & 0 & 0 &0 & 0 &0 & 0 &0 & 0 &0 & 0 &1 & 1 \\
           \color{gray}{T11} & 1 & 0 &0 & 0 &1 & 0 &0 & 0&1 & 0 &0 & 0  \\ \color{gray}{T12} & 0 & 1 &0 & 0 &0 & 1 &0 & 0 &0 & 1 &0 & 0 \\ \color{gray}{T21} & 0 & 0 &1 & 0 &0 & 0 &1 & 0 &0 & 0 &1 & 0\\ \color{gray}{T22} & 0 & 0 &0 & 1 &0 & 0 &0 & 1 &0 & 0 &0 & 1\\
\end{block}\end{blockarray}.\end{align*} 
\end{minipage}
\begin{minipage}{3cm}
\includegraphics[scale=0.55]{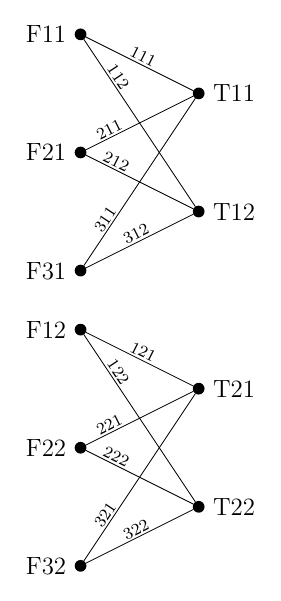}
\end{minipage}
with  circuits
\begin{align*}
\mathcal{C}_{A_{\Gamma,\mathbf r}}=\pm &\{ \scalemath{0.9}{(1,-1,0,0,-1,1,0,0,0,0,0,0), (1,-1,0,0,0,0,0,0,-1,1,0,0), (0,0,0,0,-1,1,0,0,-1,1,0,0)},\\&  \scalemath{0.9}{(0,0,1,-1,0,0,-1,1,0,0,0,0,0), (0,0,1,-1,0,0,0,0,0,0,-1,1), (0,0,0,0,0,0,1,-1,0,0,-1,1)}\}.
\end{align*}
\end{example}

\begin{figure}[!htpb]
    \centering
\scalebox{0.7}{\begin{tikzpicture}[line join = round, line cap = round] 
\coordinate [label=left:3] (3) at (0,0);
\coordinate [label=right:2] (2) at (4,0);
\coordinate [label=above:1] (1) at (2,3.5);
\begin{scope}[decoration={markings,mark=at position 0.5 with {\arrow{to}}}]
\draw[] (2)--(1)--(3)--cycle;
\end{scope}
\end{tikzpicture}
\begin{tikzpicture}[line join = round, line cap = round] 
\coordinate [label=above:3] (3) at (0,0);
\coordinate [label=right:2] (2) at (2.5,0);
\coordinate [label=right:1] (1) at (4,-2);
\coordinate [label=right:4] (4) at (4,2);
\begin{scope}[decoration={markings,mark=at position 0.5 with {\arrow{to}}}]
\draw[] (1)--(4);
\draw[fill=gray, fill opacity=.5] (2)--(3)--(4)--cycle;
\draw[fill=gray,fill opacity=.5] (2)--(1)--(3)--cycle;
\end{scope}
\end{tikzpicture}
\begin{tikzpicture}[line join = round, line cap = round] 
\coordinate [label=above:3] (3) at (0,0);
\coordinate [label=right:2] (2) at (2.5,0);
\coordinate [label=right:1] (1) at (4,-2);
\coordinate [label=right:4] (4) at (4,2);
\begin{scope}[decoration={markings,mark=at position 0.5 with {\arrow{to}}}]
\draw[fill=gray, fill opacity=.5]  (1)--(2)--(4)--cycle;
\draw[fill=gray, fill opacity=.5] (2)--(3)--(4)--cycle;
\draw[fill=gray,fill opacity=.5] (2)--(1)--(3)--cycle;
\end{scope}
\end{tikzpicture}}
    \caption{Lawrence lifts of simplicial complexes with facet sets  $\{\{1\},\{2\}\}$, $\{\{1\},\{2,3\}\}$ and $\{\{1,2\},\{2,3\}\}$, respectively. }
    \label{fig:lawrence_simplicial}
\end{figure}
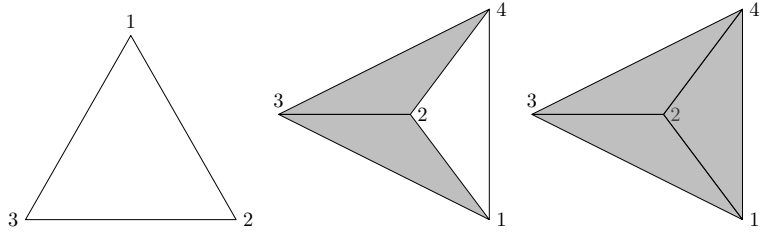

Observe that the  Lawrence lift $X_{\Lambda(A_{\Gamma,\r})}\cap \Delta_{2r-1}$ of the hierarchical model for $\Gamma,\r$  is another hierarchical model.  It is precisely  $X_{A_{\Lambda(\Gamma),(\mathbf{r},2)}}\cap \Delta_{2r-1}$  where $(\mathbf{r},2)$ is the vector $\mathbf{r}$ with entry $2$ appended to it in the last position and $\mydef{\Lambda(\Gamma)}$ is defined via its facets:  \[\mathrm{facet}(\Lambda(\Gamma))=\{[n], F\cup \{n+1\} \mid F\in \mathrm{facet}(\Gamma)\}.\] It models the dependency relations among the previous variables and an additional binary variable. \Cref{fig:lawrence_simplicial} illustrates $\Lambda(\Gamma)$ for three simplical complexes.

 The matrix $A_{\Gamma,\r}$ of a hierarchical model always contains the vector of all ones in its rowspan, since a block of rows corresponding to a fixed facet sums up to one.  Henceforth, all the circuits of~$A_{\Gamma,\r}$ are even and  \Cref{thm:mlDegMainTheoremNEW} can be applied to any~$X_{A_{\Lambda(\Gamma),(\mathbf{r},2)}}\cap \Delta_{2r-1}$  as long as the matrix $A_{\Gamma,\r}$ is totally unimodular. 
 A classification of hierarchical models with totally unimodular matrices has been completed \cite{Bernstein15, Bernstein17}. 
In the following example, we investigate some cases of hierarchical models that have a Lawrence lift structure.
\begin{example}
\label{ex:no3wayIsHier}
(i) The no-three-way interaction model in \Cref{cor:no-two-wayinteraction} is the hierarchical model\[X_{A_{(\Gamma, (m_1,m_2,2))}}\cap \Delta_{2m_1m_2-1}, \text{ where } \Gamma=2^{[3]}-[3].\] 
Hence, formulae for its degree and maximum likelihood degree are given in \Cref{cor:bipartite}. \\
(ii) The boundary of the $n$-simplex, $2^{[n]}-[n]$, is obtained by the Lawrence lift of the boundary of the $(n-1)$-simplex. Let  $\mathbf 2 \in \Z^n$ be the vector of all two-s. Then \[A_{2^{[n]}-[n], \mathbf 2}=\underbrace{\Lambda \cdots \Lambda}_{n-1\ times}(A_{K(2,2)}).\] 
The Lawrence lift preserves the property of total unimodularity, and so $A_{2^{[n]}-[n], \mathbf 2}$ is totally unimodular. Note that since $A_{K(2,2)}$ has only one circuit of  size $4$ (see Graph 2 in \Cref{ex:4graphs}), $A_{2^{[n]}-[n], \mathbf 2}$  has only one circuit of size $2^n$. We obtain,
\begin{align*}
  \deg(X_{{A}_{2^{[n]}-[n],\mathbf{2}}})&=\text{ number of spanning trees in a $2^n$-cycle }= 2^{n}, \\
  \mathrm{mldeg}(X_{{A}_{2^{[n]}-[n],\mathbf{2}}})&=\text{ number of spanning trees with external activity $0$ in a $2^n$-cycle }= 2^{n}-1. 
 \end{align*}
\end{example}

When the simplicial complex $\Gamma$ is made of exactly two facets, the matrix of the model is the incidence matrix of a bipartite graph, which we can explicitly describe and apply \Cref{thm:mld_graphs} to $\Lambda(\Gamma)$. For instance, the simplicial complex in \Cref{ex:marginalmatrix} has two facets.

\begin{proposition}
    \label{prop:hierarchical_model}
 Consider a hierarchical model with $\Gamma,\r$ where $\Gamma$ has two facets $F$ and $T$. Let  $S=F\cap T$. Denote $m_1=r_{F\setminus S}$ and $m_2=r_{T\setminus S}$.  Then,  $X_{A_{(\Gamma,\r)}}=X_{A_G}$ where  $G$ is the bipartite graph made of a disjoint union of $r_S$ complete bipartite graphs, each of the form $K(m_1, m_2)$.  
Moreover, 
\begin{align*}
  \deg(X_{A_{\Lambda(\Gamma),(\mathbf{r},2)}})&=(\tau_{K(m_1, m_2)}(1,1))^{r_S} = m_1^{r_s(m_2-1)}m_2^{r_s(m_1-1)},  \\
  \mathrm{mldeg}(X_{A_{\Lambda(\Gamma),(\mathbf{r},2)}})&=(\tau_{K(m_1, m_2)}(1,0))^{r_S} = \Big(\sum\limits_{k=1}^{m_1} \left(\frac{1}{k} \sum\limits_{i=1}^{k} (-1)^{m_1-i}{k\choose i}i^{m_1}\right)k^{m_2}\Big)^{r_s}. 
\end{align*}
\end{proposition}
\begin{proof}
Consider the bipartite graph $G$ on ground set $[\r]_{F}\sqcup [\r]_{T}$, that is, the points are labeled by tuples in $[\r]_{F}$ and $[\r]_{T}$, respectively. For each $\mathbf{i}\in [\r]_{F}$ and $\mathbf{j}\in [\r]_{T}$, the edge $(\mathbf{i},\mathbf{j})\in G$ if and only if~$\mathbf{i}|_S=\mathbf{j}|_S$. Graph $G$ has the desired properties. By  the definition of a marginal matrix we obtain  that~${A_{(\Gamma,\mathbf{r})}}=A_G$, up to an ordering of the rows.  
For the second part, note that $A_{\Lambda(\Gamma),(\mathbf{r},2)}=\Lambda(A_{\Gamma,(\mathbf{r},2)})$, and so, $X_{A_{\Lambda(\Gamma),(\mathbf{r},2)}}=X_{\Lambda(A_{G})}$. 
\Cref{cor:bipartite} concludes the argument.  
\end{proof}

\begin{example}
The matrix for  the hierarchical model in \Cref{ex:marginalmatrix} is the incidence matrix of the disjoint union of two complete graphs $K(3,2)$, as illustrated to the right of the matrix $A_{\Gamma,\mathbf r}$. Its Lawrence lift is the hierarchical model with simplical complex illustrated in \Cref{fig:lawrence_simplicial}(3) and vector of states $(3,2,2,2)$.  The Lawrence lift model has degree is $3^{2(2-1)}2^{2(3-1)}=144$ and  maximum likelihood degree is~$7^2=49$.
\end{example}

The algebraic degree of the Lawrence lift of the hierarchical model with all binary states and two disjoint facets was previously conjectured in the discussion section of \cite{coons2020generalized}, specifically, to be
\[
\deg(X_{\Lambda(A_{\Lambda(2^{[n]}\sqcup 2^{[m]}),\mathbf{2}})}) = 2^{m(2^n-1)+n(2^m-1)}.
\]
\Cref{prop:hierarchical_model} positively resolves this conjecture.

\begin{remark}
\label{rem:Gaussian}
({Relationship to Gaussian and linear models}) The M\"obius invariant appears in the literature already when discussing maximum likelihood degrees of statistical models. In \cite[Chapter 3]{sturmfels2010multivariate} the maximum likelihood degree of Gaussian models
\[
\mathcal G = \left\{K = \textbf{diag}(u_1 \cdot \x, u_2 \cdot \x, \ldots, u_n \cdot \x) \mid \x=(x_1,\ldots,x_d) \in (\mathbb{C}^\times)^d \right\} \quad \quad \text{ where } u_1,\ldots,u_n \in \mathbb{C}^d
\] with diagonal-linear concentration matrices
is determined.
This is the maximum likelihood degree of the linear model 
\[\mathbb{L} = \left\{(u_1 \cdot \x, u_2 \cdot \x, \ldots, u_n \cdot \x) \mid \x=(x_1,\ldots,x_d) \in (\mathbb{C}^\times)^d \right\} \quad \quad \text{ where } u_1,\ldots,u_n \in \mathbb{C}^d.
\]
They show that both of these maximum likelihood degrees are the \emph{M\"obius invariant} of the matroid $\mathcal M$ underlying $\mathbb{L}$. The matroid $\mathcal M$ is the column matroid of the matrix $A$ with rows $u_1,\ldots,u_n$. Thus, by \Cref{thm:MainTheorem2}, we have that the maximum likelihood degrees of three qualitatively distinct statistical models (Gaussian, linear, and log-linear) coincide:
\[
\textrm{mldeg}(X_{\Lambda(A)}) = \textrm{mldeg}(\mathcal G) = \textrm{mldeg}(\mathbb{L}) = \tau_M(1,0).
\]
We remark that the paper \cite{sturmfels2010multivariate}  erroneously calls this the \emph{Beta invariant} of $\mathcal M$. Consequently, \cite[Corollary 3.4]{sturmfels2010multivariate}  must be adapted to the condition in \Cref{prop:MLdeg1}. 
\end{remark}

\section*{Acknowledgments}
TB was partially funded by an NSERC Discovery grant (RGPIN-2023-03551).  AM was partially funded by the National Science Foundation (Grant No. DMS-2306672). The authors thank Bernd Sturmfels and Steffen Lauritzen for useful discussions in the early and late stages of the project, respectively.  
Part of this research was performed at the Max Planck Institute for Mathematics in the Sciences, University of Notre Dame, and  the Institute for Mathematical and Statistical Innovation (IMSI), which is supported by
the National Science Foundation (Grant No. DMS-1929348).

\bibliographystyle{abbrv}
\bibliography{main.bib}

\begin{thebibliography}{10}

\bibitem{LikelihoodDegenerations}
D.~Agostini, T.~Brysiewicz, C.~Fevola, L.~Kühne, B.~Sturmfels, S.~Telen, and
  T.~Lam.
\newblock Likelihood degenerations.
\newblock {\em Advances in Mathematics}, 414:108863, 2023.

\bibitem{agresti2012categorical}
A.~Agresti.
\newblock {\em Categorical data analysis}, volume 792.
\newblock John Wiley \& Sons, 2012.

\bibitem{almendra2023markov}
F.~Almendra-Hern{\'a}ndez, J.~A. De~Loera, and S.~Petrovi{\'c}.
\newblock Markov bases: a 25 year update.
\newblock {\em arXiv preprint arXiv:2306.06270}, 2023.

\bibitem{Amendola17}
C.~Am{\'e}ndola, N.~Bliss, I.~Burke, C.~Gibbons, M.~Helmer, S.~Hosten, E.~Nash,
  J.~Rodriguez, and D.~Smolkin.
\newblock The maximum likelihood degree of toric varieties.
\newblock {\em Journal of Symbolic Computation}, 92, 03 2017.

\bibitem{amendola2020maximum}
C.~Am{\'e}ndola, D.~Kosta, and K.~Kubjas.
\newblock Maximum likelihood estimation of toric fano varieties.
\newblock {\em Algebraic Statistics}, 11(1):5--30, 2020.

\bibitem{aoki2012markov}
S.~Aoki, H.~Hara, and A.~Takemura.
\newblock {\em Markov bases in algebraic statistics}, volume 199.
\newblock Springer Science \& Business Media, 2012.

\bibitem{bayer2001sysygies}
D.~Bayer, S.~Popescu, and B.~Sturmfels.
\newblock Sysygies of unimodular lawrence ideals.
\newblock 2001.

\bibitem{belay1999reye}
E.~D. Belay, J.~S. Bresee, R.~C. Holman, A.~S. Khan, A.~Shahriari, and L.~B.
  Schonberger.
\newblock Reye's syndrome in the united states from 1981 through 1997.
\newblock {\em New England Journal of Medicine}, 340(18):1377--1382, 1999.

\bibitem{Parking}
B.~Benson, D.~Chakrabarty, and P.~Tetali.
\newblock G-parking functions, acyclic orientations and spanning trees.
\newblock {\em Discrete Math.}, 310:1340--1353, 2010.

\bibitem{Bernstein15}
D.~Bernstein and S.~Sullivant.
\newblock Unimodular binary hierarchical models.
\newblock {\em Journal of Combinatorial Theory, Series B}, 123, 02 2015.

\bibitem{Bernstein17}
D.~I. Bernstein and C.~O'Neill.
\newblock Unimodular hierarchical models and their graver bases.
\newblock {\em arXiv preprint arXiv:1704.09018}, 2017.

\bibitem{bishop2007discrete}
Y.~M. Bishop, S.~E. Fienberg, and P.~W. Holland.
\newblock {\em Discrete multivariate analysis: Theory and practice}.
\newblock Springer Science \& Business Media, 2007.

\bibitem{Stiefel}
T.~Brysiewicz and F.~Gesmundo.
\newblock The degree of stiefel manifolds.
\newblock {\em Enumerative Combinatorics and Applications}, 1(3), 2019.
\newblock Article S2R20.

\bibitem{coons2020generalized}
J.~I. Coons, J.~Cummings, B.~Hollering, and A.~Maraj.
\newblock Generalized cut polytopes for binary hierarchical models.
\newblock {\em Algebraic Statistics, to appear (arXiv:2008.00043)}, 2020.

\bibitem{coons2021quasi}
J.~I. Coons and S.~Sullivant.
\newblock Quasi-independence models with rational maximum likelihood estimator.
\newblock {\em Journal of Symbolic Computation}, 104:917--941, 2021.

\bibitem{Cox11}
D.~Cox, J.~Little, and H.~Schenck.
\newblock {\em Toric Varieties}.
\newblock Graduate studies in mathematics. American Mathematical Society, 2011.

\bibitem{Dall}
A.~M. Dall.
\newblock Matroids: h-vectors, zonotopes, and lawrence polytopes.
\newblock {\em PhD diss.}, 2015.

\bibitem{de2006markov}
J.~A. De~Loera and S.~Onn.
\newblock Markov bases of three-way tables are arbitrarily complicated.
\newblock {\em Journal of Symbolic Computation}, 41(2):173--181, 2006.

\bibitem{diaconis1998algebraic}
P.~Diaconis and B.~Sturmfels.
\newblock Algebraic algorithms for sampling from conditional distributions.
\newblock {\em The Annals of statistics}, 26(1):363--397, 1998.

\bibitem{Huh}
J.~Huh.
\newblock The maximum likelihood degree of a very affine variety.
\newblock {\em Compositio Mathematica}, 149(8):1245--1266, may 2013.

\bibitem{HS14}
J.~Huh and B.~Sturmfels.
\newblock {Likelihood Geometry}.
\newblock In {\em Combinatorial algebraic geometry}, volume 2108 of {\em
  Lecture notes in mathematics}, pages 63--117. Springer, New York, 2014.

\bibitem{ishizeki2003standard}
T.~Ishizeki and H.~Imai.
\newblock Standard pairs for lawrence-type matrices and their applications to
  several lawrence-type integer programs.
\newblock {\em Optimization Methods and Software}, 18(4):417--426, 2003.

\bibitem{Kushnirenko}
A.~G. Kouchnirenko.
\newblock Poly{\'e}dres de newton et nombres de milnor.
\newblock {\em Inventiones mathematicae}, 32(1):1--31, 1976.

\bibitem{Ohsugi99}
H.~Ohsugi and T.~Hibi.
\newblock Toric ideals generated by quadratic binomials.
\newblock {\em Journal of Algebra}, 218(2):509--527, 1999.

\bibitem{oxley2006matroid}
J.~G. Oxley.
\newblock {\em Matroid theory}, volume~3.
\newblock Oxford University Press, USA, 2006.

\bibitem{Petrovic14}
S.~Petrovi{\'{c}} and D.~Stasi.
\newblock Toric algebra of hypergraphs.
\newblock {\em Journal of Algebraic Combinatorics}, 39(1):187--208, Feb 2014.

\bibitem{santos2003higher}
F.~Santos and B.~Sturmfels.
\newblock Higher lawrence configurations.
\newblock {\em Journal of Combinatorial Theory, Series A}, 103(1):151--164,
  2003.

\bibitem{ProudfootSpeyer}
D.~Speyer and N.~Proudfoot.
\newblock A broken circuit ring.
\newblock {\em Beitr{\"a}ge Algebra Geom}, 47(1):161--166, 2006.

\bibitem{Sturmfels96}
B.~Sturmfels.
\newblock {\em Grobner Bases and Convex Polytopes}.
\newblock Memoirs of the American Mathematical Society. American Mathematical
  Society, 1996.

\bibitem{sturmfels1997variation}
B.~Sturmfels and R.~R. Thomas.
\newblock Variation of cost functions in integer programming.
\newblock {\em Mathematical Programming}, 77(2):357--387, 1997.

\bibitem{sturmfels2010multivariate}
B.~Sturmfels and C.~Uhler.
\newblock Multivariate gaussians, semidefinite matrix completion, and convex
  algebraic geometry.
\newblock {\em Annals of the Institute of Statistical Mathematics},
  62(4):603--638, 2010.

\bibitem{Sullivant.2018}
S.~Sullivant.
\newblock {\em Algebraic statistics}, volume 194 of {\em Graduate Studies in
  Mathematics}.
\newblock American Mathematical Society, Providence, RI, 2018.

\bibitem{villarreal1995rees}
R.~H. Villarreal.
\newblock Rees algebras of edge ideals.
\newblock {\em Communications in Algebra}, 23(9):3513--3524, 1995.

\bibitem{vlach1986conditions}
M.~Vlach.
\newblock Conditions for the existence of solutions of the three-dimensional
  planar transportation problem.
\newblock {\em Discrete Applied Mathematics}, 13(1):61--78, 1986.

\bibitem{ZASLAVSKY198247}
T.~Zaslavsky.
\newblock Signed graphs.
\newblock {\em Discrete Applied Mathematics}, 4(1):47--74, 1982.

\end{thebibliography}

\end{document}